\documentclass[12pt,letterpaper,twoside,reqno]{amsart}
\usepackage[OT1]{fontenc}
\usepackage[left=98pt,right=98pt,top=60pt,bottom=60pt,headheight=8pt,headsep=14pt,footskip=12pt]{geometry}
\usepackage{amssymb,mathtools}
\numberwithin{equation}{section}
\usepackage{microtype}
\usepackage{booktabs,tabularx,array}
\usepackage{enumitem}
\usepackage{needspace}
\usepackage{graphicx,float}
\usepackage[hidelinks]{hyperref}
\hypersetup{pdftitle={Homothetic Copies of Polytopes and Binomial Barcodes},pdfauthor={Site Mou}}
\newtheorem{theorem}{Theorem}[section]
\newtheorem{lemma}[theorem]{Lemma}
\newtheorem{shelllemma}[theorem]{Lemma}
\newtheorem{proposition}[theorem]{Proposition}
\newtheorem{corollary}[theorem]{Corollary}
\newtheorem{conjecture}[theorem]{Conjecture}
\theoremstyle{plain}
\newtheorem{definition}[theorem]{Definition}
\newtheorem{example}[theorem]{Example}
\theoremstyle{plain}

\theoremstyle{definition}
\newtheorem*{unremark}{Remark}
\theoremstyle{plain}
\newtheorem*{introcorollary}{Corollary}
\newcommand{\R}{\mathbb R}
\newcommand{\Z}{\mathbb Z}
\newcommand{\N}{\mathbb Z_{\ge0}}
\DeclareMathOperator{\vertset}{vert}
\DeclareMathOperator{\cl}{cl}
\DeclareMathOperator{\LE}{LE}
\DeclareMathOperator{\des}{des}
\DeclareMathOperator{\conv}{conv}
\DeclareMathOperator{\vol}{vol}
\newcommand{\Ch}{\mathrm{Ch}}
\newcommand{\Ant}{\mathrm{Ant}}
\newcommand{\mset}[1]{\{\!\{#1\}\!\}}
\newcommand{\LS}[2]{\mathsf{L}_{#1}(#2)}
\setlist[enumerate]{label=(\roman*),leftmargin=*,itemsep=3pt,topsep=5pt}
\allowdisplaybreaks[2]
\title[Homothetic copies of polytopes]{Homothetic Copies of Polytopes\\and Binomial Barcodes}
\author[Site Mou]{Site Mou}
\date{}
\begin{document}
\begin{abstract}
Let \(X\subset\R^d\) be a finite full-dimensional configuration, and put
\[
F_X(N)=\max_{\substack{S\subset\R^d\\|S|=N}}
\#\{(\mathbf{v},\lambda)\in\R^d\times(0,\infty):\mathbf{v}+\lambda X\subseteq S\}.
\]
Elekes and Erd\H{o}s posed the problem of estimating \(F_X(N)\).
In this article, we determine its exact value when \(X\) is the vertex set of a simplex and obtain sharp asymptotic bounds for several classes of configurations.

We introduce binomial barcodes for order polytopes and give a recursive criterion for their effectiveness. For families in arbitrarily large dimensions, the criterion gives exact extremal values at infinitely many \(N\).
In particular, if \(X\) is the vertex set of the three-dimensional cube, then, for \(n\ge1\),
\[
F_X(n^3)=\sum_{j=1}^{n-1}j^3=\left(\frac{n(n-1)}2\right)^2.
\]
\end{abstract}
\maketitle

\section{Introduction}\label{sec:introduction}
A configuration is a finite set \(X\subset\mathbb R^d\).
We call \(X\) full-dimensional if it contains points \(\mathbf{x}_0,\ldots,\mathbf{x}_d\) such that \(\mathbf{x}_1-\mathbf{x}_0,\ldots,\mathbf{x}_d-\mathbf{x}_0\) are linearly independent. For such \(X\), write
\[
F_X(N)=
\max_{\substack{S\subset\R^d\\ |S|=N}}
\#\{(\mathbf{v},\lambda)\in\R^d\times(0,\infty):\mathbf{v}+\lambda X\subseteq S\}.
\]
For a polytope \(Q\), write \(F_Q=F_{\vertset(Q)}\); copies are counted by their vertices.

Van Kreveld and de Berg~\cite[Theorem~1]{VanKreveldDeBerg} proved that
\(F_{[0,1]^d}(N)=\Theta(N^{1+1/d})\).

In 1994, Elekes and Erd\H{o}s~\cite{ElekesErdos} considered general
full-dimensional configurations and proved
\(F_X(N)=O_X(N^{1+1/d})\).
Bernstine and Mizrahi~\cite[Theorem~3]{BernstineMizrahi} gave an
\(O(N^{1+1/d})\)-time algorithm for enumerating all positive homothetic
copies of a fixed full-dimensional configuration.

We use the functions \(K_k\) and \(\mathsf L_k\) defined by canonical
binomial expansions in \eqref{eq:canonical}--\eqref{eq:K-definition}.
For every nondegenerate \(d\)-simplex \(\Delta_d\), we prove
\(F_{\Delta_d}(N)=K_d(N)\) in Theorem~\ref{thm:simplex}.

Here \(K_d(N)\) is the Kruskal--Katona bound on \((d+1)\)-cliques in a
\(d\)-uniform hypergraph with \(N\) edges. For an order polytope, the descent
multiset of its linear extensions determines the sizes of the lattice
configurations and the number of copies they contain. Each binomial barcode
defines a function that matches these counts at the corresponding sizes
and gives a candidate upper bound for all \(N\). The recursive criterion uses inequalities for \(K_{d-1}\) to
combine bounds for the two pieces obtained by cutting along \(x_i=x_j\)
for incomparable elements \(i,j\).

For configurations admitting a regular unimodular triangulation (defined in Section~\ref{sec:asymptotic}), we
determine the optimal leading coefficient, with a linear error term.
If \(\Lambda=\operatorname{span}_{\Z}(X-\mathbf{x}_0)\), with \(\mathbf{x}_0\in X\), is a lattice of rank \(d\), write
\[
\nu=\frac{\vol(\conv X)}{\det\Lambda},
\]
where \(\det\Lambda\) is the volume of a fundamental parallelepiped of \(\Lambda\).

\Needspace{13\baselineskip}
\begin{theorem}\label{thm:unimodular}
Let \(X\subset\R^d\) be a finite full-dimensional configuration whose intrinsic lattice \(\Lambda\) has rank \(d\).
Suppose that \(P=\conv X\) has a regular unimodular triangulation with vertices in \(X\). Then
\begin{equation}\label{eq:sharp-asymptotic}
F_X(N)=\frac{1}{(d+1)\nu^{1/d}}N^{1+1/d}+O_X(N).
\end{equation}
If \(e\) denotes the number of maximal simplices in the triangulation,
then \(e=d!\nu\) and the leading coefficient is
\[
\frac{(d!)^{1/d}}{(d+1)e^{1/d}}.
\]
\end{theorem}

\begin{introcorollary}
The \(d\)-cube has \(\nu=1\), and
\[F_{[0,1]^d}(N)=\frac{1}{d+1}N^{1+1/d}+O_d(N).\]
\end{introcorollary}

For a poset \(P\) on \([d]\), its order polytope~\cite{Stanley} is
\[
O(P)=\{\mathbf{x}\in[0,1]^d:x_i\le x_j\text{ whenever }i<_Pj\}.
\]
Write \(F_P=F_{O(P)}\), and let \(e_P\) be the number of linear extensions of \(P\).
We introduce binomial barcodes \(a\) of \(P\) and associated functions
\(B_a\); see Definition~\ref{def:barcode} and
\eqref{eq:barcode-function}. Each barcode is specified by a finite
basic sequence. A barcode is effective if \(F_P(N)\le B_a(N)\) for every
\(N\ge0\).

The strict order polynomial \(\Omega_P^\circ(n)\) counts the maps
\(f:P\to\{1,\ldots,n\}\) with \(f(x)<f(y)\) whenever \(x<_P y\). Set
\(\widehat{\Omega}_P^\circ(n)=\sum_{q=0}^{n-1}\Omega_P^\circ(q)\).
Proposition~\ref{prop:lattice-endpoints} gives lattice configurations
with \(\Omega_P^\circ(n)\) points containing exactly \(\widehat{\Omega}_P^\circ(n)\) copies of \(O(P)\).
The upper bound must hold for arbitrary point sets of that size.
Comparing the coordinates of two incomparable elements splits \(O(P)\) into
two order polytopes of the same dimension whose posets have fewer linear
extensions. Theorem~\ref{thm:effective-gluing} gives a finite path test for
combining their barcode bounds.

\begin{theorem}\label{prop:endpoints}
Every barcode \(a\) of \(P\) satisfies
\[
B_a(\Omega_P^\circ(n))=\widehat{\Omega}_P^\circ(n)\qquad(n\in\N).
\]
If \(a\) is effective, then
\begin{equation}\label{eq:sharp-endpoints}
F_P(\Omega_P^\circ(n))=B_a(\Omega_P^\circ(n))=\widehat{\Omega}_P^\circ(n)\qquad(n\in\N).
\end{equation}
\end{theorem}

For the three-element antichain \(P\), we have \(\Omega_P^\circ(n)=n^3\), and
Example~\ref{ex:cube} gives an effective barcode. Theorem~\ref{prop:endpoints}
therefore yields the following exact values.
\begin{introcorollary}
For every \(n\ge1\),
\[
F_{[0,1]^3}(n^3)=\sum_{j=1}^{n-1}j^3
=\left(\frac{n(n-1)}2\right)^2.
\]
\end{introcorollary}

For the following families, the barcode bound is attained at every cardinality.
Write \(\Ch_j\) and \(\Ant_j\) for the chain and antichain on \(j\)
elements, respectively, with \(\Ch_0=\varnothing\). The symbol \(\sqcup\) denotes
disjoint union of posets; in an ordinal sum \(P\oplus Q\), every element
of \(P\) precedes every element of \(Q\).

\Needspace{14\baselineskip}
\begin{theorem}\label{thm:everywhere}
Each poset in the following table has the indicated effective basic sequence:
\begin{center}
\renewcommand{\arraystretch}{1.25}
\begin{tabular}{llll}
\toprule
Poset \(P\)&Dimension&\(e_P\)&Basic sequence\\
\midrule
\(\Ch_d\)&\(d\ge2\)&\(1\)&\((0)\)\\
\(\Ch_a\oplus\Ant_2\oplus\Ch_b\)&\(a+b+2\)&\(2\)&\((0,1)\)\\
\(\Ch_a\oplus\Ant_2\oplus\Ant_2\oplus\Ch_b\)&\(a+b+4\)&\(4\)&\((0,1,1,2)\)\\
\bottomrule
\end{tabular}
\end{center}
Here \(a\) and \(b\) are nonnegative integers. In each case the corresponding barcode \(\alpha\) satisfies
\[
F_P(N)=B_\alpha(N)\qquad(N\in\N).
\]
These equalities can be attained by nested lattice configurations.
\end{theorem}

\begin{figure}[H]
\centering
\includegraphics[width=.86\linewidth]{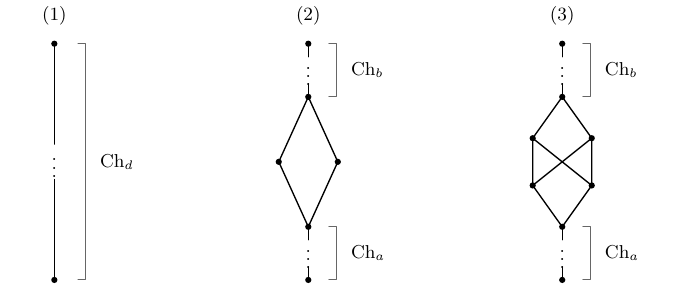}
\caption{Hasse diagrams of the families in Theorem~\ref{thm:everywhere}, in table order. Chain segments are schematic; empty chains are omitted.}
\label{fig:exact-posets}
\end{figure}

\begin{introcorollary}
For every \(N\ge1\), let \(n=\lfloor\sqrt N\rfloor\) and \(r=N-n^2\). Then
\[
F_{[0,1]^2}(N)=
\sum_{j=1}^{n-1}j^2+
\begin{cases}
\binom r2,&0\le r\le n,\\[2mm]
\binom n2+\binom{r-n}{2},&n<r<2n+1.
\end{cases}
\]
\end{introcorollary}

\begin{conjecture}\label{conj:all-effective}
Every finite poset with at least two elements admits an effective barcode.
\end{conjecture}

\begin{unremark}
Theorem~\ref{thm:effective-gluing} requires two effective barcodes and an
admissible shuffle (Definition~\ref{def:barcode-gluing}) whose path
satisfies \eqref{eq:matching-path} at every vertex. An induction based on this
theorem would require such compatible choices for a split at some incomparable
pair of each nonchain poset; see Section~\ref{sec:gluing}. Proposition~\ref{prop:enumeration} shows that the
four-dimensional cube is outside the recursive closure of chain
barcodes under this theorem. This obstruction concerns the sufficient
criterion; it does not imply that the cube is not effective.
\end{unremark}

A related problem fixes the scale and permits only translations.
In 2010, \'{A}brego, Fern\'{a}ndez-Merchant and Llano~\cite[Theorem~1]{AbregoTranslates}
obtained, for the vertex set of a nondegenerate \(d\)-simplex, the exact
maximum \(\LS{d}{N}\). After placing one pattern vertex at the origin,
their proof reduces the point set to the additive group generated by
the fixed pattern.
The reduction uses the fact that each translate lies in one coset of
that group~\cite[Section~2]{AbregoTranslates}. A positive homothetic copy
need not have this property when its scale varies, so this reduction
does not directly apply here. Theorem~\ref{thm:simplex} gives the exact
maximum \(K_d(N)\) for positive homothetic copies, using geometric layer
counting and inequalities for \(K_d\) proved in Section~\ref{sec:simplex}.

For axis-parallel squares, Kurz~\cite[p.~31, ``Point sets with many squares
determining the same directions'']{KurzSquares} recorded
computational constructions giving lower bounds at several cardinalities
up to \(17\). These bounds agree with the square formula above.
The square case of Theorem~\ref{thm:everywhere} proves their optimality
and determines the exact extremal value for every \(N\).

\vspace{\baselineskip}
Other work has studied the number of directly similar copies, allowing
rotations as well as translations and positive scaling. For a fixed
configuration \(X\subset\mathbb R^2\) with at least two points, Elekes and
Erd\H{o}s~\cite{ElekesErdos} showed that the maximum number of such copies
over all \(N\)-point subsets of \(\mathbb R^2\) is at least \(N^{2-o(1)}\).
When \(X\) has at least four points, Laczkovich and
Ruzsa~\cite{LaczkovichRuzsa} showed that this maximum is \(\Theta(N^2)\)
if and only if every complex cross-ratio of four distinct points of \(X\)
is algebraic. For squares of arbitrary orientation,
Kurz~\cite[Proposition~35]{KurzSquares} proved the upper bound
\((N^2-1)/8\).

\section{Homothetic copies of simplices}\label{sec:simplex}

We determine $F_{\Delta_d}$ for a nondegenerate $d$-simplex $\Delta_d$.

We use the canonical binomial expansions occurring in the Kruskal--Katona
theorem~\cite{Katona,Kruskal}. For each integer $d\ge1$ and $N>0$, write the unique $d$-canonical expansion as
\begin{equation}\label{eq:canonical}
N=\sum_{i=j}^{d}\binom{n_i}{i},
\qquad n_d>\cdots>n_j\ge j\ge1,
\end{equation}
and define
\begin{equation}\label{eq:K-definition}
K_d(N)=\sum_{i=j}^{d}\binom{n_i}{i+1},
\qquad \LS{d}{N}=\sum_{i=j}^{d}\binom{n_i-1}{i}.
\end{equation}
The map $\mathsf{L}_d$ is the \emph{left shift of order $d$}. Set
$K_d(0)=\LS{d}{0}=0$ and use the convention $\binom uv=0$ for $u<v$.
The canonical expansion and Pascal's identity give
\begin{equation}\label{eq:block}
\begin{gathered}
K_k\!\left(\binom mk+u\right)
=\binom m{k+1}+K_{k-1}(u),\\
\qquad 0\le u\le\binom m{k-1},\quad k\ge2,\quad m\ge k-1.
\end{gathered}
\end{equation}
They also give
\begin{equation}\label{eq:lowering}
K_d(N)=K_d(\LS{d}{N})+\LS{d}{N}.
\end{equation}
Similarly, under the hypotheses of \eqref{eq:block},
\[
\mathsf L_k\!\left(\binom mk+u\right)
=\binom{m-1}{k}+\LS{k-1}{u}.
\]
The endpoint values agree by Pascal's identity. Induction from
$\LS{1}{N}=\max\{N-1,0\}$ therefore shows that every $\mathsf L_k$
is nondecreasing.

\begin{lemma}\label{lem:superadditivity}
For all nonnegative integers $x,y$,
\[
K_d(x)+K_d(y)\le K_d(x+y).
\]
\end{lemma}

\begin{lemma}\label{lem:peak}
Let $M=\binom md$ for an integer $m\ge d$. For integers $u,v$ with
$0\le u,v\le M\le u+v$,
\[
K_d(u)+K_d(v)\le K_d(M)+K_d(u+v-M).
\]
\end{lemma}

\begin{lemma}\label{lem:tail}
Let $d\ge2$. If $\LS{d}{N}<x\le N$, then
\[
K_d(N)\ge K_d(x)+K_{d-1}(N-x).
\]
\end{lemma}

\begin{theorem}\label{thm:simplex}
For every $d\ge1$ and $N\ge0$,
\[
F_{\Delta_d}(N)=K_d(N).
\]
\end{theorem}

For the Macaulay function
\(m_d(N)=\sum_{i=j}^d\binom{n_i-1}{i-1}=N-\LS{d}{N}\),
with \(m_d(0)=0\), the inequality
\(m_{d-1}(a)+m_d(N-a)\ge m_d(N)\) holds for \(d\ge2\) and
\(0\le a<m_d(N)\) by~\cite[Lemma~2]{AbregoTranslates}.
Setting \(a=N-x\), this is equivalent to
\[
\LS{d}{N}\ge\LS{d}{x}+\LS{d-1}{N-x},
\qquad \LS{d}{N}<x\le N.
\]
Lemma~\ref{lem:tail} gives the corresponding inequality for \(K_d\).

\begin{proof}[Joint proof of Lemmas \ref{lem:superadditivity}--\ref{lem:tail} and Theorem \ref{thm:simplex}]
We induct on $d$. For $d=1$, $F_{\Delta_1}(N)=K_1(N)=\binom N2$, and
Lemmas \ref{lem:superadditivity} and \ref{lem:peak} follow directly. Assume the
assertions hold in dimensions $1,\ldots,d-1$ within their stated ranges.

We now prove, in order, Lemmas~\ref{lem:peak} and~\ref{lem:tail},
Theorem~\ref{thm:simplex}, and Lemma~\ref{lem:superadditivity} in dimension~$d$.

\medskip\noindent\textbf{Proof of Lemma \ref{lem:peak}.}

If $u=M$ or $v=M$, equality holds. If $u+v=M$, the superadditivity
of $K_d=\mathrm{old}_d$ in Frohmader~\cite[Lemma~3.2]{Frohmader}
gives $K_d(u)+K_d(v)\le K_d(M)$.
For $u,v<M<u+v$, write $u=r_d(a_d,\ldots)$, $v=r_d(c_d,\ldots)$,
and $u+v-M=r_d(b_d,\ldots)$, using Frohmader's notation for the
canonical sum in \eqref{eq:canonical}. Then
\[
r_d(a_d,\ldots)+r_d(c_d,\ldots)=r_d(m)+r_d(b_d,\ldots),
\qquad a_d,c_d<m,
\]
since $u,v<\binom md$. Apply Frohmader~\cite[Lemma~3.5]{Frohmader} with $k=d$
and $j=m$. His $r_{d+1}$ on the same binomial indices gives the
$K_d$-values in \eqref{eq:K-definition}, so
\[
K_d(u)+K_d(v)<K_d(M)+K_d(u+v-M).
\]

\medskip\noindent\textbf{Proof of Lemma \ref{lem:tail}.}

Retain the expansion \eqref{eq:canonical}, and set
\begin{equation}\label{eq:R-definition}
R=\sum_{i=j}^{d-1}\binom{n_i}{i},
\qquad \LS{d-1}{R}=\sum_{i=j}^{d-1}\binom{n_i-1}{i},
\end{equation}
with empty sums interpreted as zero. Thus
\begin{equation}\label{eq:N-split}
N=\binom{n_d}{d}+R,
\qquad \LS{d}{N}=\binom{n_d-1}{d}+\LS{d-1}{R}.
\end{equation}
Pascal's identity, the bounds $n_i\le n_d-d+i$, and the binomial
summation identity give
\begin{equation}\label{eq:R-gap}
\begin{aligned}
R-\LS{d-1}{R}&=\sum_{i=j}^{d-1}\binom{n_i-1}{i-1}\\
&\le\sum_{i=1}^{d-1}\binom{n_d-d+i-1}{i-1}
=\binom{n_d-1}{d-2}.
\end{aligned}
\end{equation}
\Needspace{8\baselineskip}
We distinguish three cases, according to whether $x\ge\binom{n_d}{d}$
and, if not, whether $N-x\le\binom{n_d-1}{d-1}$.

\noindent\textbf{Case i.} $x\ge\binom{n_d}{d}$.
By \eqref{eq:block} and Lemma \ref{lem:superadditivity} in dimension $d-1$,
applied to $x-\binom{n_d}{d}\ge0$ and $N-x\ge0$, whose sum is $R$,
\[
\begin{aligned}
K_d(N)-K_d(x)
&=K_{d-1}(R)-K_{d-1}\!\left(x-\binom{n_d}{d}\right)\\
&\ge K_{d-1}(N-x).
\end{aligned}
\]
For the remaining two possibilities, assume $x<\binom{n_d}{d}$.
Since $x>\LS{d}{N}\ge\binom{n_d-1}{d}$, applying \eqref{eq:block} to $N$ and
$x$ gives
\begin{equation}\label{eq:tail-decomposition}
K_d(N)-K_d(x)
=\binom{n_d-1}{d}+K_{d-1}(R)
-K_{d-1}\!\left(x-\binom{n_d-1}{d}\right).
\end{equation}
\noindent\textbf{Case ii.} $N-x\le\binom{n_d-1}{d-1}$.
Apply Lemma~\ref{lem:peak} in dimension $d-1$ to $N-x$ and
$x-\binom{n_d-1}{d}$ with endpoint $\binom{n_d-1}{d-1}$.
Both are nonnegative and at most this endpoint, by the case assumption
and $\binom{n_d-1}{d}\le\LS{d}{N}<x<\binom{n_d}{d}$.
Their sum is $\binom{n_d-1}{d-1}+R$ by Pascal's identity. Hence
\[
K_{d-1}(N-x)\le
\binom{n_d-1}{d}+K_{d-1}(R)
-K_{d-1}\!\left(x-\binom{n_d-1}{d}\right).
\]
This is the desired inequality by \eqref{eq:tail-decomposition}.

\noindent\textbf{Case iii.} $N-x>\binom{n_d-1}{d-1}$.
By $x>\LS{d}{N}$, \eqref{eq:N-split}, and \eqref{eq:R-gap},
\begin{equation}\label{eq:last-tail}
1\le N-x-\binom{n_d-1}{d-1}<R-\LS{d-1}{R}\le\binom{n_d-1}{d-2}.
\end{equation}
Hence $\binom{n_d-1}{d-2}>1$, which forces $d>2$. Moreover,
\eqref{eq:N-split} and $N-x>\binom{n_d-1}{d-1}$ imply
\[
\LS{d-1}{R}<x-\binom{n_d-1}{d}<R.
\]
Lemma \ref{lem:tail} in dimension $d-1\ge2$, applied to $R$ and
$x-\binom{n_d-1}{d}$ with the bounds just established, now gives
\[
\begin{aligned}
K_{d-1}(R)-K_{d-1}\!\left(x-\binom{n_d-1}{d}\right)
&\ge K_{d-2}\!\left(N-x-\binom{n_d-1}{d-1}\right).
\end{aligned}
\]
Substituting this into \eqref{eq:tail-decomposition} and using
\eqref{eq:block}, whose remainder condition is supplied by
\eqref{eq:last-tail}, yields
\[
\begin{aligned}
K_d(N)-K_d(x)
&\ge\binom{n_d-1}{d}
+K_{d-2}\!\left(N-x-\binom{n_d-1}{d-1}\right)\\
&=K_{d-1}(N-x).
\end{aligned}
\]
This proves Lemma \ref{lem:tail}.

\Needspace{4\baselineskip}
\medskip\noindent\textbf{Proof of Theorem \ref{thm:simplex}.}

For a nondegenerate simplex $\conv\{\mathbf{v}_0,\ldots,\mathbf{v}_d\}$, the affine map
\[
\mathbf{z}\longmapsto \mathbf{v}_0+\sum_{i=1}^d z_i(\mathbf{v}_i-\mathbf{v}_0)
\]
sends $\Delta_d=\conv\{\mathbf{0},\mathbf{e}_1,\ldots,\mathbf{e}_d\}$ onto it, where $\mathbf{e}_i$ are
the standard basis vectors. Thus it suffices to treat $\Delta_d$.

\emph{Upper bound.} We induct on $N$, with $N=0$ immediate. For
$S\subset\R^d$ with $|S|=N>0$, put
\[
H=\max_{\mathbf{z}\in S}\sum_{i=1}^d z_i,
\qquad x=\#\left\{\mathbf{z}\in S:\sum_{i=1}^d z_i<H\right\}.
\]
Since $x<N$, induction on $N$ bounds the copies below
$\sum_i z_i=H$ by $K_d(x)$. A copy $\mathbf{t}+\lambda\Delta_d$ meeting this
hyperplane has upper vertices $\mathbf{t}+\lambda \mathbf{e}_j$, $1\le j\le d$, and
\[
\lambda=H-\sum_{i=1}^d t_i,
\qquad t_i=\min_{1\le j\le d}(\mathbf{t}+\lambda \mathbf{e}_j)_i.
\]
It is therefore determined by its lower vertex, and also by its upper
facet. Theorem~\ref{thm:simplex} in dimension $d-1$, applied to the
$N-x$ points on the hyperplane, bounds the number of these copies by
$\min\{x,K_{d-1}(N-x)\}$.

By \eqref{eq:block}, monotonicity of $K_{k-1}$ implies monotonicity of
$K_k$ for $k\ge2$. Starting with $K_1(N)=\binom N2$, it follows that
every $K_d$ is nondecreasing.
We split according to whether $x\le\LS{d}{N}$.

\noindent\textbf{Case i.} $x\le \LS{d}{N}$.
By monotonicity and \eqref{eq:lowering},
\[
K_d(x)+\min\{x,K_{d-1}(N-x)\}
\le K_d(\LS{d}{N})+\LS{d}{N}=K_d(N).
\]
\noindent\textbf{Case ii.} $x>\LS{d}{N}$.
Lemma \ref{lem:tail} in dimension $d$, proved above, gives
\[
K_d(x)+\min\{x,K_{d-1}(N-x)\}
\le K_d(x)+K_{d-1}(N-x)\le K_d(N).
\]
The upper bound follows.

\emph{Construction attaining the bound.}
For a related lattice construction for translates of the standard simplex,
see \cite[Section~2.2]{AbregoTranslates}.
The stars-and-bars bijection~\cite[Section~1.2, pp.~18--19]{StanleyEC1}
from \(\Z_{\ge0}^d\) to the \(d\)-element subsets of \(\Z_{>0}\) is
\[
\Phi_d(\mathbf z)=\{z_1+\cdots+z_i+i:1\le i\le d\},
\]
where the \(i\)th entry is the position of the \(i\)th bar.
For distinct \(d\)-sets \(A,B\), put \(A<_{\mathrm{colex}}B\) if
\(\max(A\triangle B)\in B\), where \(\triangle\) denotes symmetric difference
\cite[Definition~2.5]{Frohmader}.
Let \(S_d(N)\) be the inverse image of the first \(N\) sets in colex order
(for their recursive description, see \cite[p.~190, proof of Theorem~1]{Katona}).
Then \(|S_d(N)|=N\) and \(S_d(N)\subseteq S_d(N+1)\).
If \(m\ge d\) and \(N\le\binom md\), the first \(N\) sets are contained
in \([m]\). Since \(\max\Phi_d(\mathbf z)=d+\sum_i z_i\), every
\(\mathbf z\in S_d(N)\) satisfies \(\sum_i z_i\le m-d\).

Order the lattice points by their images under \(\Phi_d\). This is a total
order. For a copy \(\mathbf t+\lambda\Delta_d\) with lattice vertices,
\(\lambda\) is a positive integer, and the largest vertex in this order is
\(\mathbf t+\lambda\mathbf e_1\), whose partial coordinate sums are at
least those of every other vertex. Thus each \(\mathbf z\in S_d(N)\) is
the largest vertex of exactly \(z_1\) copies: take
\(\mathbf t=\mathbf z-\lambda\mathbf e_1\) for \(1\le\lambda\le z_1\).
The other vertices precede \(\mathbf z\) and hence also lie in \(S_d(N)\).

The largest \(d\)-subset of a \((d+1)\)-set in colex is obtained by
deleting its smallest element. Since \(\min\Phi_d(\mathbf z)=z_1+1\),
the map sending the copy with largest vertex \(\mathbf z\) and scale
\(\lambda\) to \(\{\lambda\}\cup\Phi_d(\mathbf z)\) is a bijection onto
the \((d+1)\)-sets all of whose \(d\)-subsets occur among the first \(N\).
The cascade count gives exactly \(K_d(N)\) such sets
\cite[Section~3, Lemmas~3--4]{Kruskal}; see also
\cite[Section~2]{KirschRadcliffe} for the formula in colex order.
This proves Theorem~\ref{thm:simplex}.

\medskip\noindent\textbf{Proof of Lemma \ref{lem:superadditivity}.}

Take constructions attaining $K_d(x)$ and $K_d(y)$, and translate one so
that the two point sets are disjoint. Their union has $x+y$ points and
contains at least $K_d(x)+K_d(y)$ copies. Theorem \ref{thm:simplex}, just
proved in dimension $d$, gives $K_d(x)+K_d(y)\le K_d(x+y)$. This completes
the induction on $d$.

This result can also be seen in Frohmader~\cite[Lemma~3.2]{Frohmader}.
\end{proof}

\section{Regular triangulations and asymptotic bounds}\label{sec:asymptotic}
\begin{figure}[htbp]
\centering
\includegraphics[width=.84\linewidth]{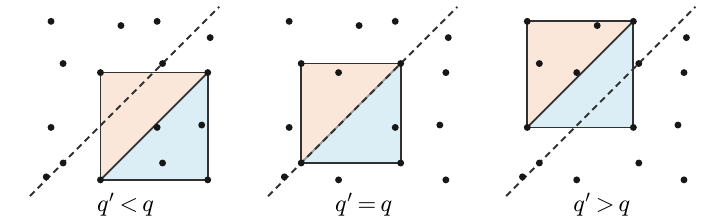}
\caption{Three squares with vertices in the same set \(S\). The dashed line is \(y-x=q=0\); each square's diagonal lies on \(y-x=q'\).}
\label{fig:cut-square}
\end{figure}

\begin{lemma}\label{lem:cut}
Let \(A,B,C\subset\R^d\) be full-dimensional convex polytopes, and let \(H\) be a hyperplane such that
\[
C=A\cup B,\qquad A\cap B=C\cap H.
\]
Suppose that \(A\) and \(B\) lie in the two respective closed half-spaces bounded by \(H\), and that
\[
\vertset(A)\cup\vertset(B)\subseteq\vertset(C).
\]
Let \(f,g:\N\to\R_{\ge0}\) be nondecreasing functions satisfying \(F_A(n)\le f(n)\) and \(F_B(n)\le g(n)\) for all \(n\ge0\). Given an \(N\)-point set \(S\) and a real number \(T\in[0,N]\), there is a hyperplane parallel to \(H\) for which the numbers \(u,v\) of points of \(S\) in the corresponding open half-spaces satisfy
\[
u\le T,\qquad v\le N-T,
\]
and the number of copies of \(C\) determined by \(S\) is at most
\begin{equation}\label{eq:cut-bound}
\min\left\{
\begin{aligned}
&f(N-v)+g(v),\\
&f(u)+g(N-u),\\
&f(u)+g(v)+\min\{u,v\}.
\end{aligned}
\right.
\end{equation}
\end{lemma}
\begin{proof}
Write \(H=\{\mathbf{x}:\boldsymbol{\eta}\cdot \mathbf{x}=h\}\), with \(A\) on the side \(\boldsymbol{\eta}\cdot \mathbf{x}\le h\). Choose \(q\in\R\) so that
\[
u:=\#\{\mathbf{s}\in S:\boldsymbol{\eta}\cdot \mathbf{s}<q\}\le T,
\qquad
v:=\#\{\mathbf{s}\in S:\boldsymbol{\eta}\cdot \mathbf{s}>q\}\le N-T.
\]
Such a value exists by ordering the numbers \(\boldsymbol{\eta}\cdot \mathbf{s}\), with repetitions retained.

For a copy \(\mathbf{t}+\lambda C\), put \(q'=\boldsymbol{\eta}\cdot \mathbf{t}+\lambda h\); its internal cutting hyperplane is \(\boldsymbol{\eta}\cdot \mathbf{x}=q'\) (Figure~\ref{fig:cut-square}). If \(q'\le q\), all vertices of \(\mathbf{t}+\lambda A\) lie in the closed lower half-space, whereas if \(q'>q\), all vertices of \(\mathbf{t}+\lambda B\) lie in the open upper half-space. Restriction to either $A$ or $B$ determines the original homothety. The vertex hypothesis therefore gives the first bound \(f(N-v)+g(v)\). Separating instead at \(q'<q\) gives the second bound \(f(u)+g(N-u)\).

For the third bound, the copies with \(q'<q\) and \(q'>q\) contribute at most \(f(u)\) and \(g(v)\), respectively. Fix vertices \(\mathbf{a}\in A\setminus H\) and \(\mathbf{b}\in B\setminus H\). If \(q'=q\), the image \(\mathbf{p}=\mathbf{t}+\lambda \mathbf{a}\) lies among the \(u\) lower points and determines the homothety, since
\[
\lambda=\frac{q-\boldsymbol{\eta}\cdot \mathbf{p}}{h-\boldsymbol{\eta}\cdot \mathbf{a}},
\qquad \mathbf{t}=\mathbf{p}-\lambda \mathbf{a}.
\]
Using \(\mathbf{b}\) in the same way gives at most \(v\) possibilities. Thus the copies with \(q'=q\) number at most \(\min\{u,v\}\), proving the third bound.
\end{proof}

We use Lemma~\ref{lem:cut} for order-polytope splits in Section~\ref{sec:gluing}.
Before that, we derive asymptotic bounds for configurations admitting
regular triangulations using coarser counting estimates.

\medskip

A triangulation of a finite full-dimensional configuration \(X\subset\mathbb R^d\) is a finite collection \(\mathcal T\) of simplices with vertices in \(X\), whose union is \(\operatorname{conv}(X)\), such that every face of a simplex in \(\mathcal T\) also belongs to \(\mathcal T\), and the intersection of any two simplices is either empty or a face of both. A simplex in \(\mathcal T\) is maximal if it is not properly contained in another simplex of \(\mathcal T\).

\begin{definition}\label{def:reg}
For a triangulation \(\mathcal T\), we call it regular if there are weights \(w:X\to\R\) and, for each maximal simplex \(\Delta\), an affine function \(\ell_\Delta\) satisfying
\[
\ell_\Delta(\mathbf{x})=w(\mathbf{x})\quad(\mathbf{x}\in\vertset(\Delta)),
\qquad
\ell_\Delta(\mathbf{x})<w(\mathbf{x})\quad(\mathbf{x}\in X\setminus\vertset(\Delta)).
\]
Equivalently, there exists a convex piecewise-affine function
\(\conv X\to\R\) whose maximal domains of affinity are precisely
the maximal simplices of \(\mathcal T\); see~\cite{Triangulations}.
\end{definition}

\begin{definition}\label{def:unimodular}
Fix \(\mathbf{x}_0\in X\) and put \(\Lambda=\operatorname{span}_{\mathbb Z}(X-\mathbf{x}_0)\). When \(\Lambda\) is a lattice of rank \(d\), it is called the intrinsic lattice of \(X\); it is independent of the choice of \(\mathbf{x}_0\).
A maximal simplex \(\Delta=\operatorname{conv}(\mathbf{v}_0,\ldots,\mathbf{v}_d)\) is unimodular with respect to \(\Lambda\) if \(\mathbf{v}_1-\mathbf{v}_0,\ldots,\mathbf{v}_d-\mathbf{v}_0\) form a \(\mathbb Z\)-basis of \(\Lambda\), equivalently if
\[
\operatorname{vol}_d(\Delta)=\frac{\det\Lambda}{d!},
\]
where \(\det\Lambda=|\det(\mathbf{b}_1,\ldots,\mathbf{b}_d)|\) for any lattice basis \(\mathbf{b}_1,\ldots,\mathbf{b}_d\). A triangulation is unimodular if all its maximal simplices are unimodular.
\end{definition}

\begin{figure}[H]
\centering
\includegraphics[width=\linewidth]{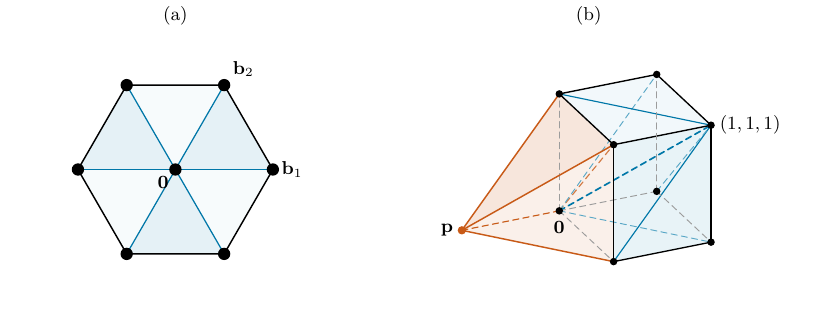}
\caption{Examples of unimodular triangulations. (a) A regular hexagon and its centre in the lattice \(\mathbb Z \mathbf{b}_1+\mathbb Z \mathbf{b}_2\), where \(\mathbf{b}_1=(1,0)\), \(\mathbf{b}_2=(1/2,\sqrt3/2)\). (b) \(\{0,1\}^3\cup\{\mathbf{p}\}\), with \(\mathbf{p}=(-1,0,0)\).}
\label{fig:unimodular-examples}
\end{figure}

\begin{theorem}\label{thm:regular}
Let \(X\subset\R^d\) be a finite full-dimensional configuration admitting a regular triangulation with \(e\) maximal simplices. Then
\begin{equation}\label{eq:regular-finite}
F_X(N)\le eK_d\!\left(\left\lfloor\frac Ne\right\rfloor\right)+\binom e2N
\qquad(N\ge0).
\end{equation}
Consequently, as \(N\to\infty\),
\[
F_X(N)\le
\frac{(d!)^{1/d}}{(d+1)e^{1/d}}N^{1+1/d}+O_X(N).
\]
\end{theorem}
\begin{proof}
Let \(N>0\), and let \(\Delta_1,\ldots,\Delta_e\) be the maximal
simplices, with \(\ell_i(\mathbf{x})=\mathbf{u}_i\cdot\mathbf{x}+c_i\)
as in Definition~\ref{def:reg}. Fix an \(N\)-point set \(S\).
Choose \(\mathbf{h}=(h_1,\ldots,h_e)\), with \(\sum_i h_i=0\), minimizing
\[
\Phi(\mathbf{h})=\sum_{\mathbf{p}\in S}\max_{1\le i\le e}
(\mathbf{u}_i\cdot\mathbf{p}+h_i).
\]
Since
\[
\Phi(\mathbf{h})\ge N\max_i h_i-
\sum_{\mathbf{p}\in S}\max_i|\mathbf{u}_i\cdot\mathbf{p}|,
\]
its sublevel sets on \(\sum_i h_i=0\) are compact, so the minimum is attained.
Put
\[
S_i=\{\mathbf{p}\in S:\mathbf{u}_i\cdot\mathbf{p}+h_i>
\mathbf{u}_j\cdot\mathbf{p}+h_j\text{ for every }j\ne i\}.
\]
For sufficiently small \(\varepsilon>0\), decrease \(h_i\) by
\((1-1/e)\varepsilon\) and increase each \(h_j\), \(j\ne i\), by
\(\varepsilon/e\). The sum of the weights remains zero, and minimality gives
\[
0\le-\left(1-\frac1e\right)\varepsilon|S_i|
+\frac{\varepsilon}{e}(N-|S_i|)
=\varepsilon\left(\frac Ne-|S_i|\right).
\]
Thus \(|S_i|\le\lfloor N/e\rfloor\).

For a copy \(\mathbf{t}+\lambda X\subseteq S\), define
\[
\theta_i=\mathbf{u}_i\cdot \mathbf{t}+h_i-\lambda c_i.
\]
Suppose first that \(\theta_i\) is the unique largest value. If \(\mathbf{x}\) is a vertex of \(\Delta_i\), then \(\ell_i(\mathbf{x})=w(\mathbf{x})\ge\ell_j(\mathbf{x})\), and hence
\[
(\mathbf{u}_i\cdot(\mathbf{t}+\lambda \mathbf{x})+h_i)
-(\mathbf{u}_j\cdot(\mathbf{t}+\lambda \mathbf{x})+h_j)
=\theta_i-\theta_j+\lambda(\ell_i(\mathbf{x})-\ell_j(\mathbf{x}))>0
\]
for every \(j\ne i\). Thus all vertices of \(\mathbf{t}+\lambda\Delta_i\) belong to \(S_i\). Theorem~\ref{thm:simplex} bounds the number of such copies by
\(K_d(|S_i|)\le K_d(\lfloor N/e\rfloor)\).

It remains to count copies for which the maximum is attained at least twice. Fix \(i<j\), and choose a vertex \(\mathbf{z}\) of \(\Delta_i\) that is not a vertex of \(\Delta_j\). Regularity gives
\(\ell_i(\mathbf{z})=w(\mathbf{z})>\ell_j(\mathbf{z})\). If \(\theta_i=\theta_j\), the point \(\mathbf{p}=\mathbf{t}+\lambda \mathbf{z}\in S\) determines both parameters, as in the proof of Lemma~\ref{lem:cut}:
\[
\lambda=\frac{(\mathbf{u}_i-\mathbf{u}_j)\cdot \mathbf{p}+h_i-h_j}
{\ell_i(\mathbf{z})-\ell_j(\mathbf{z})},
\qquad \mathbf{t}=\mathbf{p}-\lambda \mathbf{z}.
\]
There are at most \(N\) copies for this pair. Summing over the \(e\) possible unique maxima and the \(\binom e2\) pairs proves \eqref{eq:regular-finite}.

Finally, if \(a_d\) is the top index in the canonical expansion of \(n\), one has
\(n=\binom{a_d}{d}+O(a_d^{d-1})\) and
\(K_d(n)=\binom{a_d}{d+1}+O(a_d^d)\). Therefore
\[
K_d(n)=\frac{(d!)^{1/d}}{d+1} n^{1+1/d}+O(n).
\]
Substitution into \eqref{eq:regular-finite} gives the asserted asymptotic bound.
\end{proof}

\begin{proof}[Proof of Theorem~\ref{thm:unimodular}]
Translate \(X\) and \(P\) by \(-\mathbf{x}_0\), which does not change \(F_X\) or \(\nu\). Each unimodular simplex has volume \(\det\Lambda/d!\), so the triangulation has \(e=d!\nu\) maximal simplices. The upper bound follows from Theorem~\ref{thm:regular}.

For the lower bound, set \(S_m=mP\cap\Lambda\) and \(N_m=|S_m|\). Comparing \(mP\) with the fundamental cells of \(\Lambda\) gives
\[
N_m=\nu m^d+O_X(m^{d-1}).
\]
Indeed, only cells meeting a fixed-width neighbourhood of the boundary can affect the comparison, and their number is \(O_X(m^{d-1})\).

For every integer \(1\le k\le m\) and every
\(\mathbf{z}\in(m-k)P\cap\Lambda\), convexity gives
\[
\mathbf{z}+kX\subseteq(m-k)P+kP=mP.
\]
Since \(\mathbf{z}\in\Lambda\) and \(X\subseteq\Lambda\), the copy lies in \(S_m\). These are distinct homotheties, so
\[
F_X(N_m)\ge\sum_{j=0}^{m-1}N_j
=\frac{\nu}{d+1}m^{d+1}+O_X(m^d)
=\frac{N_m^{1+1/d}}{(d+1)\nu^{1/d}}+O_X(N_m).
\]
Since \(\mathbf{0}\in P\), the sets \(S_m\) are nested, and the lattice estimate also gives
\(N_{m+1}-N_m=O_X(m^{d-1})\). Given a sufficiently large \(N\), choose the largest \(m\) with \(N_m\le N\) and add arbitrary points to \(S_m\). Then \(m\asymp_X N^{1/d}\) and
\[
N^{1+1/d}-N_m^{1+1/d}
=O_X\bigl(N^{1/d}(N-N_m)\bigr)=O_X(N).
\]
Monotonicity of \(F_X\) yields the matching lower bound for every \(N\).
\end{proof}

\section{Order polytopes and binomial barcodes}\label{sec:order}
Throughout this section and the barcode arguments below, \(d\ge2\).

We use Stanley's order polytope~\cite{Stanley}. Its standard triangulation, recalled below, relates linear extensions to volume.
\begin{definition}
Let \(P\) be a poset on \(d\) elements, indexed by \([d]\). Its order polytope is
\[
O(P)=\{\mathbf{x}\in[0,1]^d:x_i\le x_j\text{ whenever }i<_Pj\}.
\]
A linear extension of \(P\) is a permutation \(\sigma\) of \([d]\) satisfying
\(\sigma^{-1}(i)<\sigma^{-1}(j)\) whenever \(i<_Pj\).
Write \(\LE(P)\) for the set of linear extensions, and put
\(e_P=|\LE(P)|\) and \(F_P=F_{O(P)}\).
\end{definition}

Choose a natural labelling, so that \(i<_Pj\) implies \(i<j\).
For a permutation \(\sigma\) of \([d]\), put
\[
\des(\sigma)=\#\{i<d:\sigma(i)>\sigma(i+1)\}.
\]
The multiset and polynomial used below are
\begin{equation}\label{eq:poset-multiset}
\mu_P=\mset{\des(\sigma):\sigma\in\LE(P)},
\qquad
Q_P(z)=\sum_{a\in\mu_P}z^a.
\end{equation}
Multiset operations and sums count multiplicities. The polynomial $Q_P$ is
the $h^*$-polynomial of $O(P)$:
\[
\sum_{m\ge0}|mO(P)\cap\Z^d|z^m
=\frac{Q_P(z)}{(1-z)^{d+1}};
\]
see \cite[Theorem~4.1]{Stanley} and
\cite[Theorem~3.15.8]{StanleyEC1}.
For a poset on a different number of elements, these definitions use its own cardinality.

\begin{proposition}\label{prop:order-triangulation}
The simplices
\[
\Delta_\sigma=
\{\mathbf{x}\in\R^d:0\le x_{\sigma(1)}\le\cdots\le x_{\sigma(d)}\le1\},
\qquad \sigma\in\LE(P),
\]
form a regular unimodular triangulation of \(O(P)\) in its intrinsic lattice \(\Z^d\). Consequently,
\[
d!\vol(O(P))=e_P,
\qquad
F_P(N)=\frac{(d!)^{1/d}}{(d+1)e_P^{1/d}}N^{1+1/d}+O_P(N).
\]
\end{proposition}
\begin{proof}
These simplices form Stanley's canonical unimodular triangulation
\cite[Section~5]{Stanley}; its regularity is proved in
\cite[Lemma~4.6]{SoprunovaSottile}.
The \(d\) edge vectors from a vertex of any maximal simplex form a basis of \(\Z^d\),
so the intrinsic lattice is \(\Z^d\). Each of the \(e_P\) maximal
simplices has volume \(1/d!\), giving the volume identity.
The asymptotic formula follows from Theorem~\ref{thm:unimodular}.
\end{proof}

\begin{proposition}\label{prop:lattice-endpoints}
For \(n\in\N\), let
\[
G_P(n)=
\{\mathbf{x}\in\{1,\ldots,n\}^d:x_i<x_j\text{ whenever }i<_Pj\}.
\]
Then
\begin{equation}\label{eq:lattice-size}
\Omega_P^\circ(n)=|G_P(n)|=\sum_{a\in\mu_P}\binom{n+a}{d}.
\end{equation}
The number of copies of \(O(P)\) whose vertices lie in \(G_P(n)\) is exactly
\begin{equation}\label{eq:lattice-copies}
\widehat{\Omega}_P^\circ(n)=\sum_{a\in\mu_P}\binom{n+a}{d+1}
=\sum_{q=0}^{n-1}\Omega_P^\circ(q).
\end{equation}
In particular, \(F_P(\Omega_P^\circ(n))\ge \widehat{\Omega}_P^\circ(n)\).
\end{proposition}
\begin{proof}
Equation~\eqref{eq:lattice-size} is the standard strict order-polynomial
expansion \cite[Theorem~3.15.8 and Corollary~3.15.12]{StanleyEC1}.

Both \(\mathbf{0}\) and \((1,\ldots,1)\) are vertices of \(O(P)\). Hence a copy \(\mathbf{t}+\lambda O(P)\) with all vertices in \(G_P(n)\) has \(\mathbf{t}\in\Z^d\), \(\lambda\in\Z_{>0}\), and \(\mathbf{t}\in G_P(n-\lambda)\). Conversely, if \(\mathbf{t}\in G_P(n-\lambda)\) and \(\lambda\) is a positive integer, then for every vertex \(\mathbf{v}\) of \(O(P)\),
\[
1\le t_i+\lambda v_i\le n,
\qquad
t_i+\lambda v_i<t_j+\lambda v_j\quad(i<_Pj),
\]
because \(v_i\le v_j\). Thus all its vertices belong to \(G_P(n)\). Distinct pairs \((\mathbf{t},\lambda)\) give distinct copies, so their number is
\[
\sum_{\lambda=1}^{n}\Omega_P^\circ(n-\lambda)=\sum_{q=0}^{n-1}\Omega_P^\circ(q).
\]
For each \(a\in\mu_P\), we have \(0\le a\le d-1\), so Pascal's identity gives
\[
\sum_{q=0}^{n-1}\binom{q+a}{d}
=\binom{n+a}{d+1}-\binom{a}{d+1}
=\binom{n+a}{d+1}.
\]
This proves \eqref{eq:lattice-copies} and the lower bound.
\end{proof}

\begin{definition}\label{def:barcode}
A barcode of \(P\) is an integer sequence \(a=(a_i)_{i\in\Z}\) satisfying
\[
a_{i+e_P}=a_i+1,\qquad
\mset{a_1,\ldots,a_{e_P}}=\mu_P.
\]
We call the finite sequence \((a_1,\ldots,a_{e_P})\) its basic sequence. The order of $a$ is $|P|-1$.
\end{definition}

For a barcode of order \(d-1\), regard each entry \(a_i\) as a bar of length \(\binom{a_i}{d-1}\) and fill the bars in order, starting with \(a_1\). A full bar $a_i$ contributes \(\binom{a_i}{d}\), and a partial filling of size \(r\) contributes \(K_{d-1}(r)\).
Thus for any positive integer $N$, write
\[
N=\sum_{1\le i<j}\binom{a_i}{d-1}+r,
\qquad
0\le r<\binom{a_j}{d-1}.
\]
We define the barcode function
\begin{equation}\label{eq:barcode-function}
B_a(N)=\sum_{1\le i<j}K_{d-1}(\binom{a_i}{d-1})+K_{d-1}(r)=\sum_{1\le i<j}\binom{a_i}{d}+K_{d-1}(r).
\end{equation}
More generally, a barcode of order \(k\ge1\) is a finite or infinite integer
sequence \(a=(a_i)\), indexed by \(\{1,\ldots,m\}\) or \(\Z_{>0}\), with bar lengths \(\binom{a_i}{k}\).
The same rule defines $B_a$ on
\[
\operatorname{dom}(B_a)=
\left\{N\in\N:N\le\sum_i\binom{a_i}{k}\right\},
\]
where the sum may be infinite. At each finite partial sum, including the
empty sum, set
\[
B_a\!\left(\sum_{i=1}^j\binom{a_i}{k}\right)
=\sum_{i=1}^j\binom{a_i}{k+1}.
\]
We call these partial sums the endpoints. Zero-length bars contribute
zero, so this convention is consistent with \eqref{eq:barcode-function}
and includes $B_a(0)=0$ and the final endpoint of a finite barcode.

Every barcode function is nondecreasing, since the partial contributions are
nondecreasing and agree with the full contributions at the endpoints.
Moreover, if $Z$ is an endpoint of an order-$k$ barcode $B$, then
\begin{equation}\label{eq:barcode-increment}
0\le B(Z+t)-B(Z)\le K_k(t)\qquad(t\in\N,\ Z+t\in\operatorname{dom}(B)).
\end{equation}
Indeed, partition the added filling into full bars and a final partial bar,
and apply Lemma~\ref{lem:superadditivity} repeatedly.

\begin{definition}\label{def:effective}
A barcode \(a\) of \(P\) is effective if
\[
F_P(N)\le B_a(N)\qquad(N\in\N).
\]
The poset \(P\), or its order polytope, is effective if it admits an effective barcode.
\end{definition}

\begin{proof}[Proof of Theorem~\ref{prop:endpoints}]
By \(a_{i+e_P}=a_i+1\), the total length of the first \(ne_P\) bars is
\[
\sum_{a\in\mu_P}\sum_{t=0}^{n-1}\binom{a+t}{d-1}
=\sum_{a\in\mu_P}\left(\binom{a+n}{d}-\binom{a}{d}\right)
=\Omega_P^\circ(n),
\]
where \(0\le a\le d-1\) makes the subtracted terms zero. Their total contribution is similarly
\[
\sum_{a\in\mu_P}\sum_{t=0}^{n-1}\binom{a+t}{d}
=\sum_{a\in\mu_P}\binom{a+n}{d+1}=\widehat{\Omega}_P^\circ(n).
\]
Zero-length bars contribute zero, so the endpoint convention in \eqref{eq:barcode-function} gives \(B_a(\Omega_P^\circ(n))=\widehat{\Omega}_P^\circ(n)\). For effective \(a\), this identity and Proposition~\ref{prop:lattice-endpoints} give~\eqref{eq:sharp-endpoints}.
\end{proof}

We next compare barcode functions under adjacent swaps, using the left
shift \(\LS{k}{r}\) defined in \eqref{eq:K-definition}.

\begin{proposition}\label{prop:left-transfer}
Let $m\ge k\ge1$ be integers and put $L=\binom{m}{k}$.
If integers $r,x\ge0$ satisfy $\LS{k}{r}<x\le L$, then
\begin{equation}\label{eq:left-transfer}
K_k(x)+K_k(r+L-x)\le K_k(L)+K_k(r).
\end{equation}
\end{proposition}
\begin{proof}
We distinguish whether $r\le x$.

\noindent\textbf{Case i.} $r\le x$. Then $x,r+L-x\le L$ and their sum is
$L+r\ge L$. Apply Lemma~\ref{lem:peak} to these two arguments with endpoint $L$.
\noindent\textbf{Case ii.} $r>x$. Then $k\ge2$, since $\LS{1}{r}=r-1$.
Write $r=\sum_{i=s}^k\binom{n_i}{i}$ in $k$-canonical form.
The inequality $\LS{k}{r}<L$ implies $n_k\le m$. Hence
\[
r-\LS{k}{r}=\sum_{i=s}^k\binom{n_i-1}{i-1}
\le\sum_{i=1}^k\binom{m-k+i-1}{i-1}=\binom{m}{k-1}.
\]
It follows that $0<r-x<\binom{m}{k-1}$.
Equation \eqref{eq:block} and Lemma~\ref{lem:tail} in
dimension $k$, applied with $N=r$, now give
\[
K_k(x)+K_k(r+L-x)
=K_k(L)+K_k(x)+K_{k-1}(r-x)
\le K_k(L)+K_k(r).
\]
\end{proof}

For finite integer sequences $u,v$ of the same length, write
$v\preccurlyeq u$ if $v$ is obtained from $u$ by finitely many adjacent
swaps $(x,y)\mapsto(y,x)$ with $x>y$, allowing zero swaps.
For sequences indexed by $\Z_{>0}$, write $v\preccurlyeq u$ if, for every
$m\ge1$, finitely many such swaps in $u$ produce a sequence whose first
$m$ entries are $v_1,\ldots,v_m$.
For barcodes of $P$, compare their positive-indexed parts; then
\[
a\preccurlyeq a'
\quad\Longleftrightarrow\quad
(a_1,\ldots,a_{e_P})\preccurlyeq(a'_1,\ldots,a'_{e_P}).
\]

\begin{lemma}\label{lem:reordering}
Let $a,a'$ be barcodes of the same order $k\ge1$. If $a\preccurlyeq a'$, then
\[
 B_a(N)\le B_{a'}(N)\qquad(N\in\operatorname{dom}(B_a)).
\]
\end{lemma}
\begin{proof}
First consider one allowed swap. For its labels $x>y$, put $X=\binom xk$ and
$Y=\binom yk$, so $X\ge Y$. If $X=Y$, the functions agree.
Let $T$ be the endpoint preceding the two bars. The functions agree for
$N\le T$ and $N\ge T+X+Y$. Otherwise, put $s=N-T$.
We split $0\le s\le X+Y$ at $Y$ and $X$.

\noindent\textbf{Case i.} $0\le s\le Y$. Both contributions are $K_k(s)$.

\noindent\textbf{Case ii.} $Y\le s\le X$. Lemma~\ref{lem:superadditivity} gives
\[
 K_k(Y)+K_k(s-Y)\le K_k(s).
\]
\noindent\textbf{Case iii.} $X\le s\le X+Y$. Since $Y,s-Y\le X$ and their
sum is $s\ge X$, Lemma~\ref{lem:peak} with endpoint $X$ gives
\[
 K_k(Y)+K_k(s-Y)\le K_k(X)+K_k(s-X).
\]
All other contributions agree.
For general $a\preccurlyeq a'$, fix $N\in\operatorname{dom}(B_a)$ and choose
$m$ so that the first $m$ bars of $a$ have total length at least $N$.
Finitely many allowed swaps in $a'$ give $a_1,\ldots,a_m$ as its first $m$
entries. Applying the preceding inequalities at each swap gives
$B_a(N)\le B_{a'}(N)$.
\end{proof}

In particular, if $a,a'$ are barcodes of $P$ with $a\preccurlyeq a'$,
then $a'$ is effective whenever $a$ is.

\section{Poset decompositions and barcode shuffles}\label{sec:gluing}

Let $C$ be a poset on a set $S$ of cardinality $d\ge2$, and let $x,y\in S$
be incomparable. Following \cite[Sections~1--2]{StanleyDecomposition}, define
\[
 A=\cl(C\cup\{x<y\}),\qquad
 B=\cl(C\cup\{y<x\}),\qquad D=C/(x=y),
\]
where $\cl$ denotes transitive closure. More precisely, let
$\pi:S\to(S\setminus\{x,y\})\cup\{*\}$ identify $x$ and $y$ with $*$.
The strict order of $D$ is the transitive closure of
$\{(\pi(u),\pi(v)):u<_C v\}$.
A minimal directed cycle would pass through $*$ and lift to a strict
chain in $C$ with both endpoints in $\{x,y\}$, contradicting
irreflexivity or the incomparability of $x,y$. Thus $D$
is a poset on $d-1$ elements. We call $(A,B,D)$ the
\emph{order-polytope split} of $C$ at $x,y$: geometrically,
\[
O(C)=O(A)\cup O(B),\qquad O(A)\cap O(B)=O(C)\cap\{t_x=t_y\}.
\]
Interchange $x$ and $y$,
and hence $A$ and $B$, when necessary so that $e_A\le e_B$.

\begin{proposition}\label{prop:root-identity}
An order-polytope split satisfies $e_C=e_A+e_B$ and $e_D\le e_A\le e_B$. Moreover,
\begin{equation}\label{eq:root-polynomial}
 Q_C(z)=Q_A(z)+Q_B(z)+(z-1)Q_D(z).
\end{equation}
Here $Q_D$ is defined using the dimension $d-1$ of $O(D)$.
\end{proposition}
\begin{proof}
The extensions of $C$ with $x$ before $y$ are exactly $\LE(A)$, and those
with $y$ before $x$ are exactly $\LE(B)$. Thus $e_C=e_A+e_B$.
Replacing $*$ in the order specified by an extension of $D$ by the consecutive
pair $x,y$ gives an extension of $A$: every relation of $C$ is respected,
as is $x<y$. This map is injective. Replacing $*$ by $y,x$ similarly
injects $\LE(D)$ into $\LE(B)$. Hence $e_D\le\min\{e_A,e_B\}$.

The strict order-polynomial recurrence
\cite[Section~2]{StanleyDecomposition} gives
\[
\Omega_C^\circ(n)=\Omega_A^\circ(n)+\Omega_B^\circ(n)+\Omega_D^\circ(n).
\]
The strict order-polynomial generating function
\cite[Theorem~3.15.8 and Corollary~3.15.12]{StanleyEC1} is
\[
\sum_{n\ge0}\Omega_P^\circ(n)z^n
=\frac{z^{|P|}Q_P(z^{-1})}{(1-z)^{|P|+1}}.
\]
Since $A,B,C$ have $d$ elements and $D$ has $d-1$, substitution gives
\[
Q_C(z^{-1})=Q_A(z^{-1})+Q_B(z^{-1})
+(z^{-1}-1)Q_D(z^{-1}),
\]
which is \eqref{eq:root-polynomial}.
\end{proof}

\begin{definition}\label{def:barcode-gluing}
Let $a,b$ be barcodes of $A,B$. An $(e_A,e_B)$-shuffle is a permutation
$\pi$ of $[e_C]$ satisfying
\[
\pi(1)<\cdots<\pi(e_A),\qquad
\pi(e_A+1)<\cdots<\pi(e_C).
\]
It defines a sequence $c$ by
\begin{equation}\label{eq:shifted-blocks}
\begin{aligned}
c_{\pi(i)}&=a_{e_D+i} &&(1\le i\le e_A),\\
c_{\pi(e_A+j)}&=b_j &&(1\le j\le e_B),
\end{aligned}
\end{equation}
extended by $c_{h+e_C}=c_h+1$.
The shuffle $\pi$ is \emph{admissible} if $c$ is a barcode of $C$.
\end{definition}

The path of $\pi$ starts at $(0,0)$, takes horizontal steps in positions
$\pi(1),\ldots,\pi(e_A)$ and vertical steps in the remaining positions,
and ends at $(e_A,e_B)$.

By \eqref{eq:root-polynomial} and $a_{i+e_A}=a_i+1$,
\[
\sum_{h=1}^{e_C}z^{c_h}-Q_C(z)
=(z-1)\left(\sum_{i=1}^{e_D}z^{a_i}-Q_D(z)\right).
\]
Thus $\pi$ is admissible if and only if
\begin{equation}\label{eq:gluing-prefix}
\mset{a_1,\ldots,a_{e_D}}=\mu_D.
\end{equation}

\begin{definition}\label{def:matching}
Let $e_A\ge2$ and $m=e_A-1$. For $k,M\in\Z$, the entries $a_k$ and $b_M$ are
\emph{matched} if there exist $R\ge M+m-1$ and a finite integer sequence $w$ such that
\[
 (b_M,\ldots,b_R)\preccurlyeq(a_{k-1},\ldots,a_{k-m},w).
\]
Here $w$ denotes the sequence appended after $a_{k-m}$, rather than a
single entry. For $e_A=1$, every pair $(a_k,b_M)$ is matched.
\end{definition}

In Figure~\ref{fig:cube-gluing}(b), both barcodes have basic sequence
\((0,1,1)\), and \(e_D=2\). At the starting vertex \((0,0)\), \(a_3\)
and \(b_1\) are matched:
\((b_1,b_2)=(0,1)\preccurlyeq(1,0)=(a_2,a_1)\).

\Needspace{13\baselineskip}
\begin{unremark}
To test whether $a_k$ and $b_M$ are matched when $e_A>1$, put
\[
J=\min\{j\ge M:b_j>\max_{1\le i<e_A}a_{k-i}\}.
\]
The relation $b_{j+e_B}=b_j+1$ ensures that $J$ exists.
In Definition~\ref{def:matching} it suffices to test $R<J$.
If $R\ge J$, the occurrence $b_J$ cannot pass any of $a_{k-1},\ldots,a_{k-m}$
under the allowed swaps. Deleting the occurrences ending at positions
$J,\ldots,R$ from both sequences gives the comparison with $R=J-1$.
For each $R<J$, the entries $b_M,\ldots,b_R$ fix the multiset of $w$, and
testing $\preccurlyeq$ is a finite search over permutations.

Starting with chain barcodes, enumerate splits, admissible shuffles and
reorderings. Each poset has finitely many barcodes and $e_A,e_B<e_C$
for each split, so the search decides whether an effective barcode for $C$
can be certified by Theorem~\ref{thm:effective-gluing} and
Lemma~\ref{lem:reordering}.
\end{unremark}

\begin{theorem}\label{thm:effective-gluing}
Let $(A,B,D)$ be an order-polytope split of $C$, with $e_A\le e_B$.
Let $a,b$ be effective barcodes of $A,B$, and let $\pi$ be an admissible
shuffle with associated barcode $c$.
Suppose that
\begin{equation}\label{eq:matching-path}
a_{e_D+i+1}\text{ and }b_{j+1}\text{ are matched}
\end{equation}
at every vertex $(i,j)$ of the path of $\pi$.
Then $c$ is effective:
\[
F_C(N)\le B_c(N)\qquad(N\in\N).
\]
\end{theorem}
\begin{proof}
Extend the path of \(\pi\) by successive translates by \((e_A,e_B)\).
If \(e_A\ge2\), a comparison in Definition~\ref{def:matching} at a vertex \((i,j)\)
gives one at \((i+e_A,j+e_B)\) by increasing every entry of both sequences
by one, since \(a_{s+e_A}=a_s+1\) and \(b_{s+e_B}=b_s+1\).
The replacement \((u,v)\mapsto(v,u)\) with \(u>v\) becomes
\((u+1,v+1)\mapsto(v+1,u+1)\), with \(u+1>v+1\).
Thus \eqref{eq:matching-path} holds on the extended path. For \(e_A=1\)
it holds by definition.

For \(N\in\N\), choose \(h\ge0\) and \(r\) so that
\[
N=\sum_{s=1}^{h}\binom{c_s}{d-1}+r,
\qquad 0\le r<\binom{c_{h+1}}{d-1}.
\]
Let \((i,j)\) be the vertex after \(h\) steps, and put
\[
k=e_D+i,\qquad
U_t=\sum_{s\le t}\binom{a_s}{d-1},\qquad
V=\sum_{s=1}^{j}\binom{b_s}{d-1}.
\]
By \eqref{eq:shifted-blocks}, \eqref{eq:gluing-prefix}, and \eqref{eq:barcode-function},
\[
N=U_k+V+r,\qquad
B_c(N)=B_a(U_k)+B_b(V)+K_{d-1}(r).
\]

We now apply Lemma~\ref{lem:cut} to \(O(A),O(B),O(C)\).
If \(x,y\) are the incomparable elements defining the split, then
\[
O(A)=O(C)\cap\{\mathbf t:t_x\le t_y\},\qquad
O(B)=O(C)\cap\{\mathbf t:t_y\le t_x\}.
\]
These order polytopes are full-dimensional, and
\(O(A)\cap O(B)=O(C)\cap H\) for \(H=\{\mathbf t\in\R^d:t_x=t_y\}\).
By \cite[Corollary~1.3]{Stanley}, \(\vertset(O(P))=O(P)\cap\{0,1\}^d\)
for \(P=A,B,C\). Hence
\[
\vertset(O(A))\cup\vertset(O(B))
=(O(A)\cup O(B))\cap\{0,1\}^d
=\vertset(O(C)).
\]
Take \(f=B_a\), \(g=B_b\), and \(T=U_k+r=N-V\).
The functions \(B_a,B_b\) are nondecreasing, and effectiveness gives
\(F_A\le B_a\) and \(F_B\le B_b\).
For an arbitrary \(N\)-point set, the lemma supplies a hyperplane
parallel to \(H\) with \(u\le U_k+r\) and \(v\le V\) points in the
open half-spaces corresponding to \(O(A)\) and \(O(B)\), respectively.
By \eqref{eq:cut-bound}, it
suffices to prove
\[
\min\left\{
\begin{aligned}
&B_a(N-v)+B_b(v),\\
&B_a(u)+B_b(N-u),\\
&B_a(u)+B_b(v)+\min\{u,v\}
\end{aligned}
\right\}
\le B_c(N).
\]

We use the third term when \(u\le U_{k-e_A}+\mathsf L_{d-1}(r)\)
and the second term otherwise.

\noindent\textbf{Case i.} \(u\le U_{k-e_A}+\mathsf L_{d-1}(r)\).
Shifting the barcode indices by \(e_A\) and using Pascal's identity gives
\[
B_a(U_k)-B_a(U_{k-e_A})
=\sum_{s\le k-e_A}
\left(\binom{a_s+1}{d}-\binom{a_s}{d}\right)
=U_{k-e_A}.
\]
By monotonicity of \(B_a\), \eqref{eq:barcode-increment}, and
\eqref{eq:lowering},
\[
\begin{aligned}
B_a(u)+u
&\le B_a(U_{k-e_A})+U_{k-e_A}
   +K_{d-1}\!\bigl(\mathsf L_{d-1}(r)\bigr)+\mathsf L_{d-1}(r)\\
&=B_a(U_k)+K_{d-1}(r).
\end{aligned}
\]
Since \(v\le V\), the third term is at most
\(B_a(u)+B_b(V)+u\le B_c(N)\).

\Needspace{5\baselineskip}
Assume now that \(u>U_{k-e_A}+\mathsf L_{d-1}(r)\).
To bound the second term, we first use the matching condition.

If \(e_A\ge2\), take \(R\) and \(w\) from Definition~\ref{def:matching}
for the matched pair \((a_{k+1},b_{j+1})\). The adjacent transpositions
in the finite comparison act only on the first \(R-j\) entries. Appending
\(b_{R+1},b_{R+2},\ldots\) to both sequences therefore gives
\[
(b_{j+1},b_{j+2},\ldots)
\preccurlyeq
\widetilde b:=(a_k,\ldots,a_{k-e_A+2},w,b_{R+1},b_{R+2},\ldots).
\]
If \(e_A=1\), take \(\widetilde b=(b_{j+1},b_{j+2},\ldots)\).
For \(k-e_A+1\le t\le k\), the first \(k-t\) entries of \(\widetilde b\)
give the endpoint \(U_k-U_t\), with
\(B_{\widetilde b}(U_k-U_t)=B_a(U_k)-B_a(U_t)\).

The barcode function of \((b_{j+1},b_{j+2},\ldots)\) is
\(s\mapsto B_b(V+s)-B_b(V)\). Lemma~\ref{lem:reordering}
and \eqref{eq:barcode-increment} at \(U_k-U_t\) give, for
\(k-e_A+1\le t\le k\) and \(z\in\mathbb Z_{\ge0}\),
\[
\begin{aligned}
B_b(V+U_k-U_t+z)-B_b(V)
&\le B_{\widetilde b}(U_k-U_t+z)\\
&\le B_a(U_k)-B_a(U_t)+K_{d-1}(z).
\end{aligned}
\tag{*}
\]

Let \(t\) be the largest index in \(\{k-e_A,\ldots,k\}\) with
\(U_t\le u\), and set \(s=u-U_t\).
Then \(B_a(u)\le B_a(U_t)+K_{d-1}(s)\) by \eqref{eq:barcode-increment}.

\noindent\textbf{Case ii.} \(t>k-e_A\) and \(s\le r\).
Taking \(z=r-s\) in \((*)\) at index \(t\) gives
\[
B_a(u)+B_b(N-u)
\le B_a(U_k)+B_b(V)+K_{d-1}(s)+K_{d-1}(r-s)
\le B_c(N),
\]
where the last inequality is Lemma~\ref{lem:superadditivity}.

\noindent\textbf{Case iii.} \(t=k-e_A\) or \(s>r\).
Then \(t<k\), since \(t=k\) would give \(s=u-U_k\le r\).
Moreover,
\[
\mathsf L_{d-1}(r)<s<\binom{a_{t+1}}{d-1}.
\]
The lower bound follows from the hypothesis on \(u\) if \(t=k-e_A\)
and from \(s>r\ge\mathsf L_{d-1}(r)\) otherwise; the upper bound
follows from \(u<U_{t+1}\).
Take \(z=U_{t+1}-u+r\) in \((*)\) at index \(t+1\).
Since \(B_a(U_{t+1})-B_a(U_t)=\binom{a_{t+1}}d\), this gives
\[
\begin{aligned}
&B_a(u)+B_b(N-u)-B_a(U_k)-B_b(V)\\
&\qquad\le K_{d-1}(s)+K_{d-1}\!\left(\binom{a_{t+1}}{d-1}+r-s\right)
-\binom{a_{t+1}}d
\le K_{d-1}(r),
\end{aligned}
\]
by Proposition~\ref{prop:left-transfer} with \(m=a_{t+1}\).
Thus \(B_a(u)+B_b(N-u)\le B_c(N)\).
\end{proof}

\Needspace{11\baselineskip}
\section{Exact bounds for order polytopes}\label{sec:examples}

\begin{shelllemma}\label{lem:shell-filling}
Let $Q\subset\R^d$, $d\ge2$, be a full-dimensional lattice polytope with a
unimodular simplex facet $F$. For $\mathbf{v}\in\Z^d$ and $s\in\N$, put
$B=(\mathbf{v}+sF)\cap\Z^d$. Let $S$ be a finite set disjoint from $B$.
Suppose that whenever $\mathbf{t}+\lambda F$ has all its vertices in $B$, with
$\lambda>0$, every remaining vertex of $\mathbf{t}+\lambda Q$ belongs to $S$.
Then the points of $B$ can be ordered so that, for every integer
$0\le r\le |B|$, the first $r$ points together with $S$ determine at least
$K_{d-1}(r)$ copies of $Q$ not already contained in $S$.
\end{shelllemma}
\begin{proof}
Choose a vertex $\mathbf{p}$ of $F$. With origin $\mathbf{v}+s\mathbf{p}$
and the edge vectors at $\mathbf{p}$ as a lattice basis, $B$ has coordinates
$\{\mathbf{z}\in\Z_{\ge0}^{d-1}:\sum_i z_i\le s\}$. Under $\Phi_{d-1}$,
these points correspond to all the $(d-1)$-subsets of $[s+d-1]$, a colex
initial segment. In the induced order, the first $r$ points correspond
to $S_{d-1}(r)$ in the proof of Theorem~\ref{thm:simplex}; at
$r=|B|=\binom{s+d-1}{d-1}$, they fill $B$. Each of the $K_{d-1}(r)$ facet
copies determines a unique copy of $Q$. These copies are distinct and
were not contained in $S$, since their facet vertices lie in $B$.
\end{proof}

\begin{proof}[Proof of Theorem~\ref{thm:everywhere}]
\emph{Effectiveness.} We treat the three families in the table in order.

\noindent\textbf{Case i.} For $P=\Ch_d$, the basic sequence $(0)$ is effective
by Theorem~\ref{thm:simplex}.

\noindent\textbf{Case ii.} Let $P=\Ch_a\oplus\Ant_2\oplus\Ch_b$, and split
at the antichain pair. Then $A\cong B\cong\Ch_d$, $D\cong\Ch_{d-1}$,
and $e_A=e_B=e_D=1$. By Case~i, $A$ and $B$ have effective barcodes with
basic sequence $(0)$. The first entry for $A$ gives $\mset{0}=\mu_D$,
so \eqref{eq:gluing-prefix} holds. Since $e_A=1$, \eqref{eq:matching-path} holds by
Definition~\ref{def:matching}. The shuffle in Figure~\ref{fig:family-gluing}(a)
gives an effective barcode with basic sequence $(0,1)$ by
Theorem~\ref{thm:effective-gluing}.

\noindent\textbf{Case iii.} Let
$P=\Ch_a\oplus\Ant_2\oplus\Ant_2\oplus\Ch_b$, and split at the first
antichain pair. Then $A\cong B\cong\Ch_{a+2}\oplus\Ant_2\oplus\Ch_b$
and $D\cong\Ch_{a+1}\oplus\Ant_2\oplus\Ch_b$, with $e_A=e_B=e_D=2$.
By Case~ii, $A$ and $B$ have effective barcodes with basic sequence $(0,1)$.
The first two entries for $A$ have multiset $\mset{0,1}=\mu_D$,
so \eqref{eq:gluing-prefix} holds.
Along the path in Figure~\ref{fig:family-gluing}(b), the comparisons in
Definition~\ref{def:matching} are $(0,1)\preccurlyeq(1,0)$ at $(0,0)$,
$(1)\preccurlyeq(1)$ at $(0,1),(1,1),(1,2)$, and
$(1,2)\preccurlyeq(2,1)$ at $(2,2)$. Theorem~\ref{thm:effective-gluing}
gives an effective barcode with basic sequence $(0,1,1,2)$.

\begin{figure}[H]
\centering
\includegraphics[width=\linewidth]{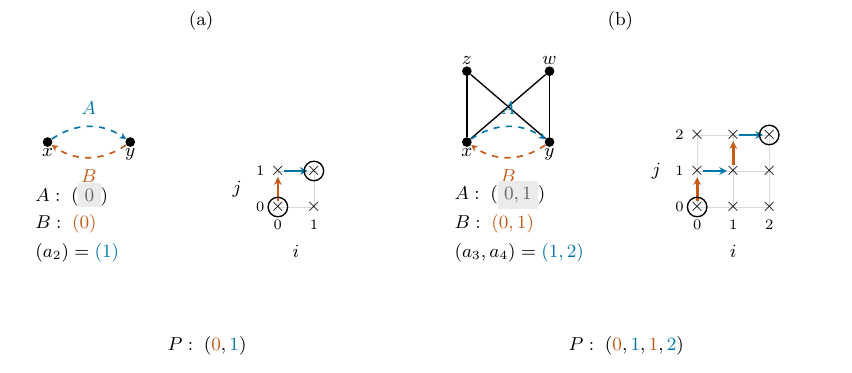}
\caption{Splits and shuffle paths for Theorem~\ref{thm:everywhere}.
A blue step to $(i,j)$ uses $a_{e_D+i}$; an orange step uses $b_j$.
The first $e_D$ entries of $a$ are gray. Crosses mark vertices satisfying
\eqref{eq:matching-path}; dashed arrows show the added relations.}
\label{fig:family-gluing}
\end{figure}

\smallskip\noindent\emph{Attaining configurations.}
The chain case follows from Theorem~\ref{thm:simplex}.
Put $c=a+b$. For $b>0$, the integral affine bijection in ordinal-sum coordinates
\[
\begin{gathered}
(\xi_1,\ldots,\xi_{d-b},y_1,\ldots,y_b)\\
\longmapsto
(1-y_b,\ldots,1-y_1,\xi_1+1-y_1,\ldots,\xi_{d-b}+1-y_1)
\end{gathered}
\]
sends $O(P)$ onto $O(\Ch_c\oplus\Ant_2)$ or
$O(\Ch_c\oplus\Ant_2\oplus\Ant_2)$, respectively, without changing $F_P$.
For $b=0$, use the identity. Write $\mathbf{e}_1,\ldots,\mathbf{e}_d$
for the standard basis and $\mathbf{x}=(x_1,\ldots,x_c)$.

For the polytope $Q$ in each case, set $G_m=mQ\cap\Z^d$ for $m\ge0$.
We partition each layer $G_{m+1}\setminus G_m$ into ordered blocks.
To apply Lemma~\ref{lem:shell-filling}, we show that each copy of $Q$
with its specified facet vertices in the current block has all its
other vertices in $G_m$ or in earlier blocks of that layer.

\smallskip\noindent\textbf{Case $e_P=2$.}
Take $Q=O(\Ch_c\oplus\Ant_2)$:
\[
Q=\left\{(\mathbf{x},u_1,u_2):
0\le x_1\le\cdots\le x_c\le\min\{u_1,u_2\},\quad u_1,u_2\le1\right\}.
\]
An empty chain is omitted, retaining $0\le\min\{u_1,u_2\}$.
The facets $F_i=Q\cap\{u_i=1\}$ are unimodular simplices.
The blocks of $G_{m+1}\setminus G_m$, in order, are:
\[
\begin{array}{c|c}
\text{condition}&\text{block}\\ \hline
u_1>u_2&(\mathbf{e}_{c+1}+mF_1)\cap\Z^d\\
u_1\le u_2&(m+1)F_2\cap\Z^d
\end{array}
\]

Suppose $\mathbf{t}+\lambda F_i$, $\lambda>0$, has its vertices in
block~$i$. The vertices $\mathbf{e}_{c+i}$ and
$\mathbf{e}_{c+1}+\mathbf{e}_{c+2}$ of $F_i$ give
$\lambda\in\Z_{>0}$, $\mathbf{t}\in\Z^d$, and
\[
0\le t_1\le\cdots\le t_c\le t_{c+3-i}\le t_{c+i}=m+1-\lambda.
\]
Thus $\vertset(\mathbf{t}+\lambda Q)\subseteq G_{m+1}$.
For $\mathbf{z}\in\vertset(Q)\setminus F_i$,
\[
(\mathbf{t}+\lambda\mathbf{z})_{c+i}=m+1-\lambda\le m.
\]
In block~1, $t_{c+2}+\lambda\le m$, so all remaining vertices lie
in $G_m$. In block~2, any remaining vertex outside $G_m$ has
$u_1=m+1>u_2$ and lies in block~1.

\smallskip\noindent\textbf{Case $e_P=4$.}
The affine map
\[
(\mathbf{x},u_1,u_2,w_1,w_2)
\longmapsto(\mathbf{x},u_1,u_2,1-w_1,1-w_2)
\]
sends $O(\Ch_c\oplus\Ant_2\oplus\Ant_2)$, with chain coordinates
$\mathbf{x}$ and antichain pairs $(u_1,u_2),(w_1,w_2)$, onto
\[
Q=\left\{(\mathbf{x},u_1,u_2,v_1,v_2):
\begin{array}{l}
0\le x_1\le\cdots\le x_c\le\min\{u_1,u_2\},\\
v_1,v_2\ge0,\quad u_i+v_j\le1\quad(i,j=1,2)
\end{array}
\right\}.
\]
On $F_{ij}=Q\cap\{u_i+v_j=1\}$, eliminating $v_j$ gives
\[
0\le x_1\le\cdots\le x_c\le u_{3-i}\le u_i\le1-v_{3-j}\le1,
\]
so $F_{ij}$ is a unimodular simplex.
The blocks of $G_{m+1}\setminus G_m$, in order, are as follows; omit row~1
when $m=0$:
\[
\begin{array}{c|c}
\text{conditions}&\text{block}\\ \hline
u_1>u_2,\ v_1>v_2&(\mathbf{e}_{c+1}+\mathbf{e}_{c+3}+(m-1)F_{11})\cap\Z^d\\
u_1>u_2,\ v_1\le v_2&(\mathbf{e}_{c+1}+mF_{12})\cap\Z^d\\
u_1\le u_2,\ v_1<v_2&(\mathbf{e}_{c+4}+mF_{22})\cap\Z^d\\
u_1\le u_2,\ v_1\ge v_2&(m+1)F_{21}\cap\Z^d
\end{array}
\]

Suppose $\mathbf{t}+\lambda F_{ij}$, $\lambda>0$, has its vertices
in the current block. The vertices $\mathbf{e}_{c+i}$ and
$\mathbf{e}_{c+1}+\mathbf{e}_{c+2}$ give
$\lambda\in\Z_{>0}$ and $\mathbf{t}\in\Z^d$. Put
$\mathbf{q}=\mathbf{t}+\lambda(\mathbf{e}_{c+1}+\mathbf{e}_{c+2})$.
The facet inequalities imply
\[
0\le t_1\le\cdots\le t_c\le t_{c+3-i}\le t_{c+i},
\qquad t_{c+3},t_{c+4}\ge0,
\]
and, for $r,s=1,2$,
\[
\max_{\mathbf{z}\in\mathbf{t}+\lambda Q}(z_{c+r}+z_{c+2+s})
=t_{c+r}+t_{c+2+s}+\lambda
=q_{c+r}+q_{c+2+s}\le m+1.
\]
Hence $\vertset(\mathbf{t}+\lambda Q)\subseteq G_{m+1}$.
The row conditions at $\mathbf{q}$ therefore give, on $\mathbf{t}+\lambda Q$,
\[
\begin{array}{c|l}
F_{11}&\max\{u_1+v_2,u_2+v_1,u_2+v_2\}\le m,\\
F_{12}&\max\{u_2+v_1,u_2+v_2\}\le m,\\
F_{22}&u_2+v_1\le m.
\end{array}
\]
These inequalities exclude every later block. For
$\mathbf{z}\in\vertset(Q)\setminus F_{ij}$, its $0,1$ coordinates give
$z_{c+i}+z_{c+2+j}=0$, so
\[
(\mathbf{t}+\lambda\mathbf{z})_{c+i}
+(\mathbf{t}+\lambda\mathbf{z})_{c+2+j}
=m+1-\lambda\le m.
\]
Thus each remaining vertex lies in $G_m$ or an earlier block.

Apply Lemma~\ref{lem:shell-filling} successively to these blocks.
A block of scale $s$ has $\binom{s+d-1}{d-1}$ points and bar label
$s+d-1$. The successive labels, including zero-length bars, are
\[
\begin{array}{c|l}
e_P=2&m+d-1,\ m+d,\\
e_P=4&m+d-2,\ m+d-1,\ m+d-1,\ m+d.
\end{array}
\]
The initial set $G_0=\{\mathbf{0}\}$ supplies a bar of label $d-1$.
The basic sequences are therefore $(0,1)$ and
$(0,1,1,2)$, respectively. Their effectiveness gives equality at
every cardinality.

\end{proof}

\Needspace{7\baselineskip}
Now we consider the square and the three-dimensional cube.

\begin{example}\label{ex:square}
The square has basic sequence \((0,1)\). For \(N\ge1\), write \(N=n^2+r\), where \(n\ge1\) and \(0\le r<2n+1\). Then
\[
F_{[0,1]^2}(N)=
\sum_{j=1}^{n-1}j^2+
\begin{cases}
\binom r2,&0\le r\le n,\\[2mm]
\binom n2+\binom{r-n}{2},&n<r<2n+1.
\end{cases}
\]
\end{example}
\begin{proof}
Apply Theorem~\ref{thm:everywhere} to \(P=\Ant_2\), using \(K_1(t)=\binom t2\).
\end{proof}

\begin{figure}[H]
\centering
\includegraphics[width=\linewidth]{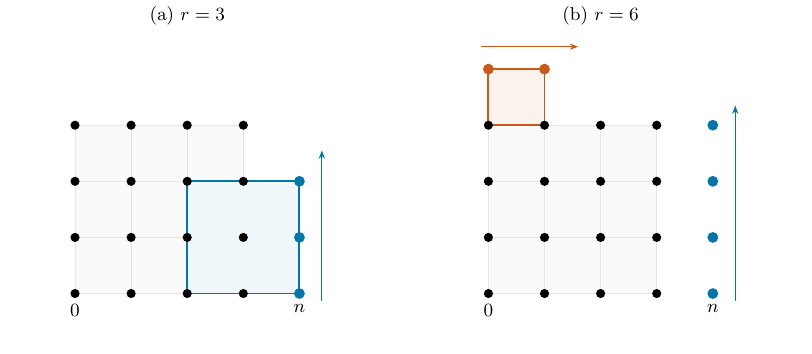}
\caption{Attaining configurations for \(n=4\): (a) \(r=3\), (b) \(r=6\). One new square is highlighted in each panel.}
\label{fig:square-extremal}
\end{figure}

\begin{unremark}
Figure~\ref{fig:square-extremal} illustrates the attaining configurations.
Start with $\{0,\ldots,n-1\}^2$ and insert
$(n,0),\ldots,(n,n-1)$, followed by $(0,n),\ldots,(n,n)$.
Each pair in the added column, and then each pair in the added row,
determines a new square whose other two vertices have already been
inserted. The increments are therefore $\binom r2$ for $r\le n$, and
$\binom n2+\binom{r-n}{2}$ for $r>n$.
\end{unremark}

\Needspace{8\baselineskip}
\begin{example}\label{ex:cube}
The basic sequence $(1,0,1,1,2,1)$ is effective for $[0,1]^3$.
Consequently,
\[
F_{[0,1]^3}(n^3)=\sum_{j=1}^{n-1}j^3
=\left(\frac{n(n-1)}2\right)^2\qquad(n\ge1).
\]
The $n\times n\times n$ grid attains these values.
\end{example}
\begin{proof}
\noindent\textbf{(1)} Split $\Ch_2\sqcup\Ch_1$ at the isolated element
and the smaller chain element. Then $A\cong\Ch_3$ and
$B\cong\Ch_1\oplus\Ant_2$, with common facet identified with $O(\Ch_2)$.
Take effective barcodes $a,b$ of $A,B$ with basic sequences $(0)$ and $(0,1)$,
respectively, by Theorem~\ref{thm:everywhere}. We have
$\mset{a_1}=\mu_{\Ch_2}=\mset{0}$, and $e_A=1$ gives \eqref{eq:matching-path}.
The shuffle in Figure~\ref{fig:cube-gluing}(a) gives an effective barcode
with basic sequence $(b_1,a_2,b_2)=(0,1,1)$ by Theorem~\ref{thm:effective-gluing}.

\noindent\textbf{(2)} Split $\Ant_3$ at two elements. Then
$A\cong B\cong\Ch_2\sqcup\Ch_1$, with common facet identified with
$O(\Ant_2)$. Take effective barcodes $a,b$ of $A,B$ with basic sequence
$(0,1,1)$ from (1). We have $\mset{a_1,a_2}=\mu_{\Ant_2}=\mset{0,1}$.
The match at \((0,0)\) was checked after Definition~\ref{def:matching}.
On the path in Figure~\ref{fig:cube-gluing}(b), \eqref{eq:matching-path}
at \((1,0),(3,2),(3,3)\) follows, respectively, from the adjacent swaps
\[
(1,1,0)\longmapsto(0,1,1),\qquad
(2,1,1)\longmapsto(1,1,2),\qquad
(2,1)\longmapsto(1,2).
\]
At $(1,1),(2,1),(2,2)$, the comparison in Definition~\ref{def:matching} is
$(1,1)\preccurlyeq(1,1)$. Theorem~\ref{thm:effective-gluing} thus gives an
effective barcode with basic sequence
\[
(a_3,b_1,a_4,b_2,a_5,b_3)=(1,0,1,1,2,1).
\]
\begin{figure}[H]
\centering
\includegraphics[width=\linewidth]{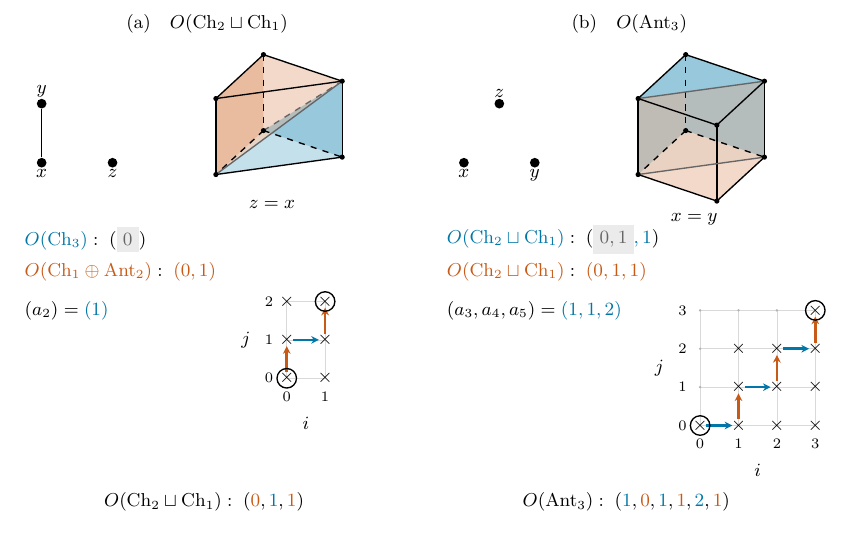}
\caption{Splits and shuffle paths in Example~\ref{ex:cube}.
Blue entries are from $a$; orange entries are from $b$.}
\label{fig:cube-gluing}
\end{figure}

Since $G_{\Ant_3}(n)=\{1,\ldots,n\}^3$ and $\Omega_{\Ant_3}^\circ(n)=n^3$,
Theorem~\ref{prop:endpoints} gives
$F_{[0,1]^3}(n^3)=\widehat{\Omega}_{\Ant_3}^\circ(n)=\sum_{q=0}^{n-1}q^3$.
\end{proof}

\begin{unremark}
Every three-dimensional order polytope is effective. Up to isomorphism,
the three-element posets are $\Ch_3$, $\Ch_1\oplus\Ant_2$,
$\Ant_2\oplus\Ch_1$, $\Ch_2\sqcup\Ch_1$ and $\Ant_3$.
The first three cases follow from Theorem~\ref{thm:everywhere} and the
remaining two from Example~\ref{ex:cube}.
\end{unremark}

\Needspace{17\baselineskip}
In dimensions $4,5,6$, we enumerate the order polytopes certified effective
by recursive applications
of Theorem~\ref{thm:effective-gluing} and Lemma~\ref{lem:reordering}, starting with
chain barcodes.

\begin{proposition}\label{prop:enumeration}
The exhaustive enumeration of this recursive closure gives the following counts.
\begin{center}
\renewcommand{\arraystretch}{1.2}
\begin{tabular}{ccc}
\toprule
Dimension&All congruence classes&Certified classes\\
\midrule
\(4\)&\(12\)&\(11\)\\
\(5\)&\(39\)&\(31\)\\
\(6\)&\(183\)&\(97\)\\
\bottomrule
\end{tabular}
\end{center}
The last column counts congruence classes certified effective by these recursive rules.
\end{proposition}
\begin{proof}
The supplementary program performs the following enumeration and verification.

\smallskip\noindent
(1) Generate all posets on $d=4,5,6$ elements by successively adding a
maximal element above each order ideal and identifying isomorphic posets.
Their upper-set indicator vectors are the vertices of their order
polytopes. Group these vertex sets by bijections preserving squared
Euclidean distances.

\smallskip\noindent
(2) Starting from the chain sequence $(0)$, process posets by increasing
$e_P$. For each split, take $e_A\le e_B$, using both orders of $A,B$ if
$e_A=e_B$. Consider all admissible shuffles of certified barcodes of $A,B$
whose paths satisfy
\eqref{eq:matching-path}. Theorem~\ref{thm:effective-gluing} certifies
each resulting barcode. Include the reorderings supplied by
Lemma~\ref{lem:reordering}.
Since $e_A<e_P$ and $e_B<e_P$, this computes the recursive
closure specified above.

\smallskip\noindent
(3) Recheck each stored certificate by verifying the certificates for $A,B$,
the split $(A,B,D)$, the reorderings, and \eqref{eq:matching-path} at
every path vertex.
Count a congruence class as certified if it contains a poset with a
verified barcode. This gives the table.
\end{proof}

\Needspace{12\baselineskip}
\begin{corollary}\label{cor:enumerated-extrema}
For each class certified in Proposition~\ref{prop:enumeration}, the
recursive enumeration gives exact values at the cardinalities in
Theorem~\ref{prop:endpoints}. For example, for the posets \(P\) and \(Q\)
shown below and every \(n\in\N\),
\begin{center}
\begin{tabular}{@{}c@{\qquad}l@{}}
\begin{tabular}{@{}c@{}}
\includegraphics{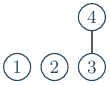}\\[1pt]
\(P\) (\(d=4,\ e_P=12\))
\end{tabular}
&
\(\displaystyle
F_P\!\left(\frac{n^3(n-1)}2\right)
=\frac{n(n-1)(n-2)(12n^2-9n-1)}{120}
\)
\\
\addlinespace[8pt]
\begin{tabular}{@{}c@{}}
\includegraphics{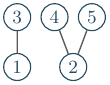}\\[1pt]
\(Q\) (\(d=5,\ e_Q=20\))
\end{tabular}
&
\(\displaystyle
F_Q\!\left(\frac{n^2(n-1)^2(2n-1)}{12}\right)
=\frac{n^2(n-1)^2(n-2)^2}{36}\,.
\)
\end{tabular}
\end{center}
\end{corollary}

\Needspace{6\baselineskip}
When even the shorter period \(e_A\) is long, a match may fail because it
requires entries from the preceding period of \(a\). For instance, if
\(e_A\ge e_D+i+2\), the block of \(e_A-1\) entries at \((i,j)\) reaches
into that period and includes the zero from \eqref{eq:gluing-prefix}.
A \(1\) shifted down by \(a_{s-e_A}=a_s-1\) would give a second zero,
but \(b_1,b_2,\ldots\) has only one. In that case no vertex in column \(i\)
is matched, and every shuffle path for these barcodes meets that column.

\begin{figure}[H]
\centering
\includegraphics[width=\linewidth]{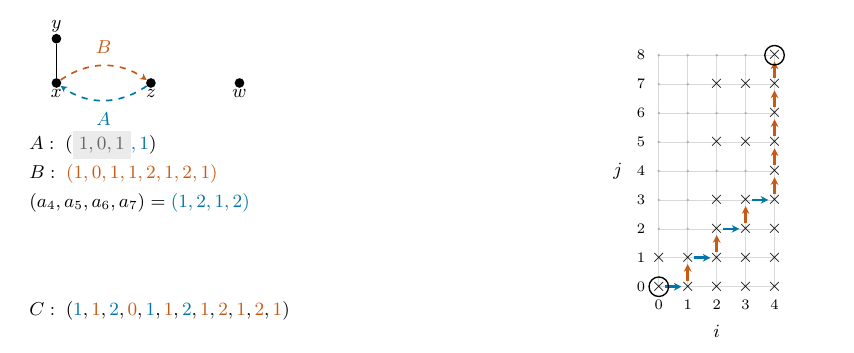}
\caption{A shuffle path certifying effectiveness for \(\Ch_2\sqcup\Ant_2\).}
\label{fig:four-dimensional-success}
\end{figure}

\begin{figure}[H]
\centering
\includegraphics[width=\linewidth]{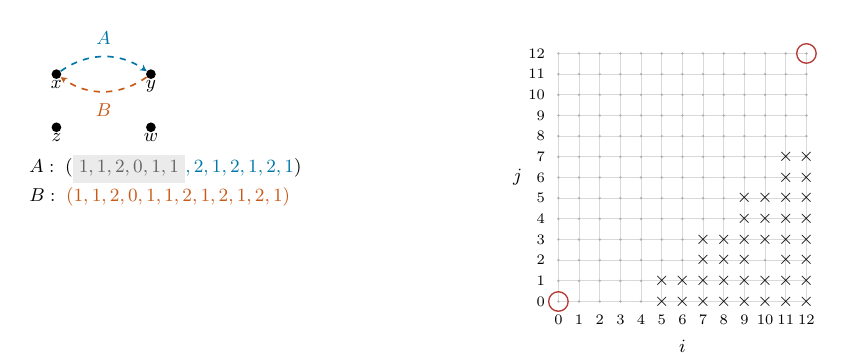}
\caption{The matching grid for the displayed barcode pair of \(\Ant_4\).}
\label{fig:four-dimensional-failure}
\end{figure}

We conclude with a family in arbitrary dimension.

\begin{example}\label{ex:unbounded}
For \(d\ge2\) and \(1\le m\le d\), the poset
\[
P_{d,m}=(\Ch_{m-1}\sqcup\{z\})\oplus\Ch_{d-m}
\]
admits an effective barcode with basic sequence
\[
(0,\underbrace{1,\ldots,1}_{m-1}).
\]
Moreover,
\[
e_{P_{d,m}}=m,\qquad \vol(O(P_{d,m}))=\frac{m}{d!}.
\]
These \(d\) order polytopes are pairwise noncongruent. Thus, if \(M_d\)
counts the \(d\)-dimensional congruence classes certified by the recursive
criterion, then \(M_d\ge d\).
\end{example}

\begin{figure}[H]
\centering
\includegraphics[width=.50\linewidth]{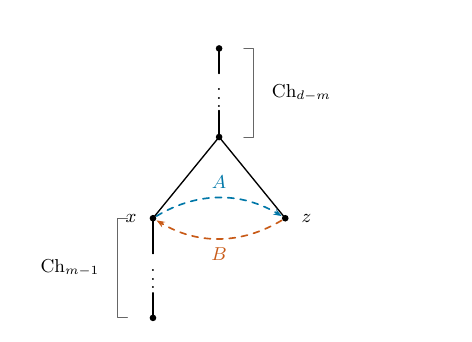}
\caption{The schematic split of \(P_{d,m}\), \(m\ge2\), at \(x=\max(\Ch_{m-1})\) and \(z\); empty chains are omitted.}
\label{fig:unbounded-poset}
\end{figure}
\begin{proof}
Induct on \(m\). For \(m=1\), the poset is a chain and has sequence
\((0)\). For \(m\ge2\), take the order-polytope split of \(P_{d,m}\)
at \(x=\max(\Ch_{m-1})\) and \(z\), with \(A\) adding
\(x<z\) and \(B\) adding \(z<x\). Then
\(A\cong\Ch_d\), \(B\cong P_{d,m-1}\), and
\(D\cong\Ch_{d-1}\). Let \(a\) be the chain barcode of \(A\)
and \(b\) the inductively obtained barcode of \(B\). Since
\(e_A=e_D=1\), we have \(\mset{a_1}=\mset{0}=\mu_D\), and
\eqref{eq:matching-path} holds by Definition~\ref{def:matching}.
Place \(a_2\) after \(b_1,\ldots,b_{m-1}\).
Theorem~\ref{thm:effective-gluing} gives an effective barcode with basic
sequence \((0,1,\ldots,1)\).

A linear extension is determined by the position of \(z\) among the \(m-1\) elements of the first chain, so \(e_{P_{d,m}}=m\). Proposition~\ref{prop:order-triangulation} gives volume \(m/d!\). Distinct values of \(m\) give distinct volumes and hence noncongruent polytopes, proving \(M_d\ge d\).
\end{proof}

In other words, Theorem~\ref{prop:endpoints} gives, for every \(n\in\N\),
\[
F_{P_{d,m}}\!\left(\binom{n}{d}+(m-1)\binom{n+1}{d}\right)
=\binom{n}{d+1}+(m-1)\binom{n+1}{d+1}.
\]

\medskip\noindent\textbf{Supplementary material.}
The accompanying file \path{Homothetic_polytopes_computation.zip}
contains \texttt{compute.py}, \texttt{catalog.json},
\texttt{certificates.jsonl.gz} and \texttt{README.html}.
With Python~3.8 or later, running \texttt{compute.py} and entering a
dimension from $2$ to $6$ regenerates and verifies the corresponding
results. The HTML file provides an offline explorer of the saved posets,
barcode representatives and certificates for Theorem~\ref{thm:effective-gluing}, and is updated after
each successful run.

\section*{Declaration of generative AI and AI-assisted technologies}
During the research and preparation of this article, the author used
ChatGPT (OpenAI) to assist with finite enumeration and computational checks, development and optimization
of the supplementary code, preparation of figure-generation code, and
the organization and language of the manuscript. The figures were drawn
from explicit mathematical descriptions using Python and TikZ.
The mathematical ideas and proofs were developed by the author,
who reviewed and revised the AI-assisted material and takes full
responsibility for the article and supplementary material.

\end{document}